\documentclass[preprint,12pt]{elsarticle}

\usepackage{amssymb,latexsym,amsmath}
\usepackage{amsthm}

\renewcommand{\bullet}{\cdot}

\newcommand{\Z}{\mathsf{Z}}
\newcommand{\disp}{\displaystyle}
\newcommand{\dist}{{\rm dist}}
\newcommand{\wid}{\widetilde}
\newcommand{\vg}{\vadjust{\goodbreak}}
\newenvironment{pf}{\vskip-\lastskip\medskip
\noindent{\bf Proof.} }{\vskip-\lastskip\vskip4pt plus2pt}

\renewenvironment{itemize}[1]{\begin{list}{}%
{\def\makelabel##1{\hss\llap{##1}}%
\setlength{\topsep}{4pt}%
\setlength{\parsep}{0pt}%
\setlength{\itemsep}{0pt}%
\settowidth{\labelwidth}{#1\enspace}%
\setlength{\leftmargin}{\labelwidth}
}}
{\end{list}}

\makeatletter
\def\ps@pprintTitle{%
     \let\@oddhead\@empty
     \let\@evenhead\@empty
     \def\@oddfoot{\footnotesize\itshape
\hfill\today}%
     \let\@evenfoot\@oddfoot}

\def\eqalign#1{\null\,\vcenter{\openup\jot\m@th
 \ialign{\strut\hfil$\displaystyle{##}$&$\displaystyle{{}##}$\hfil
 \crcr#1\crcr}}\,}

\makeatother

\font\sy=bbold10 at14truept
\newcommand{\1}{\hbox{\sy 1}}

\newcommand{\tekstl}[2]{\vskip-\lastskip\vskip4pt{\indent
\begin{minipage}[b]{11.9cm}#1\end{minipage}}\hfill
{\rm\refstepcounter{equation}(\theequation)\label{#2}}%
\vskip-\lastskip\vskip4pt}

\def\<#1>{\langle#1\rangle}

\renewcommand{\Delta}{\varDelta}
\renewcommand{\Phi}{\varPhi}

\newcommand{\hRe}{\mathop{\rm Re}\nolimits}

\newcommand{\calL}{\mathcal L}
\newcommand{\calB}{\mathcal B}
\newcommand{\calC}{\mathcal C}
\newcommand{\calD}{\mathcal D}
\newcommand{\calE}{\mathcal E}
\newcommand{\calF}{\mathcal F}
\newcommand{\calM}{\mathcal M}
\newcommand{\calN}{\mathcal N}
\newcommand{\calO}{\mathcal O}
\newcommand{\calR}{\mathcal R}
\newcommand{\calS}{\mathcal S}

\newcommand{\bbC}{\mathbb C}
\newcommand{\bbN}{\mathbb N}
\newcommand{\bbR}{\mathbb R}

\newcommand{\bbZ}{\mathbb Z}

\renewcommand{\[}{\mathclose[}
\renewcommand{\]}{\mathopen]}

\renewcommand{\theequation}{\thesection.\arabic{equation}}
\makeatletter
\@addtoreset{equation}{section}
\makeatother

\begin{document}
 
\parskip0pt
\parindent10pt

\makeatletter

\def\ps@copyright{}

\def\th@plain{%
 \thm@preskip4pt plus2pt minus2pt
 \thm@postskip4pt plus2pt minus2pt
 \itshape 
 }

\def\printFirstPageNotes{%
  \iflongmktitle
	  \let\columnwidth=\textwidth\fi
  \ifx\@tnotes\@empty\else\@tnotes\fi
  \ifx\@cornotes\@empty\else\@cornotes\fi
  \ifx\@elseads\@empty\relax\else
   \let\thefootnote\relax
   \footnotetext{\textit{Email address:\space}%
     \@elseads}\fi
  \ifx\@elsuads\@empty\relax\else
   \let\thefootnote\relax
   \footnotetext{\textit{URL:\space}%
     \@elsuads}\fi
  \ifx\@fnotes\@empty\else\@fnotes\fi
  \iflongmktitle\if@twocolumn
	  \let\columnwidth=\Columnwidth\fi\fi
} 
 \makeatother

\theoremstyle{plain}
\newtheorem{theor}{Theorem}
\newtheorem{coro}{Corollary}
\newtheorem{propo}{Proposition}
\newtheorem{lem}{Lemma}

\newtheorem{Def}{Definition}
\newtheorem{rem}{Remark}
\newtheorem{exa}{Example}

\begin{frontmatter}

\title{The Petrovski\u\i\  correctness and semigroups of operators}
\author{Jan Kisy\'nski}
\address{Institute of Mathematics, Polish Academy of Sciences,\\
\'Sniadeckich 8, 00-956 Warszawa, Poland}
\ead{jan.kisynski@gmail.com} 

\begin{abstract}
Let $P(\partial/\partial x)$ be an $m\times n$ matrix whose entries are PDO on $\bbR^n$ with constant coefficients, and let $\calS(\bbR^n)$ be the space of infinitely differentiable rapidly decreasing functions on $\bbR^n$. It is proved that $P(\partial/\partial x)|_{(\calS(\bbR^n))^m}$ is the infinitesimal generator of a $(C_0)$-semigroup
$(S_t)_{t\ge0}\subset L((\calS(\bbR^n))^m)$ if and only if $P(\partial/\partial x)$  satisfies the Petrovski\u\i\  correctness condition. Moreover, if it is the case, then $(S_t)_{t\ge0}$ is an exponential semigroup whose characteristic exponent is equal to the stability index of $P(\partial/\partial x)$. Similar statements are also proved for some other function spaces on $\bbR^n$, and for the space of tempered distributions.

\bigskip\noindent
{\it MSC:} Primary 35E15, 47D06; Secondary 46F05, 15A42
 
\end{abstract}

\begin{keyword}Cauchy's problem; Petrovski\u\i\  correct system; $(C_0)$-semigroup 
\end{keyword}

\end{frontmatter}

\section{Introduction.
The matricial differential operator  $P(\partial/\partial x)$  and   the
    corresponding Cauchy problem}\label{sec1}

Let $m,n,d\in\bbN$. Denote by $M_m$ the ring of $m\times m$ matrices with complex entries. Suppose that for every multiindex $\alpha=(\alpha_1,\ldots,\alpha_n)\in\bbN^n_0$ of length $|\alpha|= \alpha_1+\cdots
+\alpha_n\le d$
there is given a matrix $A_{\alpha}\in M_{m}$. Consider the polynomial
$P(X)=P(X_1,\ldots,X_n)\in M_m[X_1,\ldots,X_n]$
of  $n$  variables $X_1,\ldots,X_n$ with coefficients in $M_m$ defined by the formula 
\begin{equation}\label{eq1.1}
P(X)=\sum_{|\alpha|\le d}A_{\alpha} X^{\alpha},\ \quad
X^\alpha=X^{\alpha_1}_1\cdots X^{\alpha_n}_n.
\end{equation}
Substituting
$$
X_k=\frac{\partial}{\partial x_k},\ \quad k=1,\ldots,n,
$$
we obtain the matricial partial differential operator on $\bbR^n$
with constant coefficients:
\begin{equation}\label{eq1.2}
P\bigg(\frac{\partial}{\partial x}\bigg)
=\sum_{|\alpha|\le d}A_\alpha\bigg(\frac{\partial}{\partial x}\bigg)^\alpha,\ \quad
\bigg(\frac{\partial}{\partial x}\bigg)^\alpha=\bigg(\frac{\partial}{\partial x_1}\bigg)^{\alpha_1}\cdots\bigg(\frac{\partial}{\partial x_n}\bigg)^{\alpha_n}.
\end{equation}
Substituting
$$
X_k=i\xi_k\in\bbR,\ \quad \xi_k\in\bbR,\ \quad k=1,\ldots,n,
$$
we obtain the symbol of $P(\partial/\partial x)$, i.e. the 
$m\times m$ matrix
\begin{equation}\label{eq1.3}
\wid{P}=\wid{P}(\xi)=\sum_{|\alpha|\le d}i^{|\alpha|} A_{\alpha}\xi^{\alpha},\ \quad
\xi^{\alpha}=\xi^{\alpha_1}_{1}\cdots\xi^{\alpha_n}_{n}, 
\end{equation}
whose entries are scalar complex polynomials on~$\bbR^n$.

If  $E$  is a space of $C^m$-valued functions or distributions on~$\bbR^n$,   then one
can consider the Cauchy problem for  $E$-valued functions $u(\bullet)$ of a real
variable:
\begin{equation}\label{eq1.4}
\eqalign{
\frac{du(t)}{dt}&= P\bigg(\frac{\partial}{\partial x}\bigg) u(t)\ \quad
\mbox{for }t\ge0\ \mbox{(or for $t\in\bbR$)}, \cr
u(0)&=u_0.\cr}
\end{equation}
where $u_0\in E$ is given. Some sort of well posedness of such a Cauchy problem consists in the fact that the operator $P(\partial/\partial x)$ considered on the domain
$\{u\in E:P(\partial/\partial x)u\in E\}$
is the infinitesimal generator of a one-parameter semigroup
$(S_t)_{t\ge0}\subset L(E)$ (or group $(G_t)_{t\in\bbR}\subset L(E)$) of class $(C_0)$.

The subsequent Section~\ref{sec2} is devoted to $E=(Z'_n)^m$. Then the space of the
Fourier transforms $\calF^{-1}E=(\calF^{-1}Z'_n)^m=(\calD'(\bbR^n))^m$ is invariant with respect to multiplication by arbitrary elements of $C^\infty(\bbR^n;M_m)$,
and this implies that the operator $P(\partial/\partial x)|_{(Z'_n)^m}$
is the infinitesimal generator of a one-parameter group 
$(G_t)_{t\in\bbR}\subset L((Z'_n)^m)$ of class~$(C_0)$. The main result of the present paper  is formulated in Section~\ref{sec3} where several spaces  $E$  are considered with 
$\calF^{-1}E$   not invariant with respect to multiplication by arbitrary elements 
of~$C^\infty(\bbR^n;M_m)$. Then, in order to prove that a suitable restriction of
$P(\partial/\partial x)$ generates  a $(C_0)$-semigroup 
$(S_t)_{t\ge 0}\subset L(E)$, one must assume {\it something} about $\wid P$, and, for each of the spaces  $E$  considered, this {\it something} appears to be the Petrovski\u\i\  correctness condition.

Recall that if  $E$  is an  l.c.v.s.~and $L(E)$ is the algebra of continuous linear operators on~$E$, then a parametrized family $(S_t)_{t\ge 0}\subset L(E)$ is called a
{\it one-parameter semigroup of class} $(C_0)$ if it satisfies the following three conditions:
\begin{itemize}{\rm(iii)}
\item[\rm(i)] $S_{t_1+t_2}=S_{t_1}S_{t_2}$ for every $t_1,t_2\in[0,\infty\[$,
\item[\rm(ii)] $S_0=\1$, the unity of $L(E)$,
\item[\rm(iii)]  for every $u\in E$ the map $[0,\infty\[\ni t\mapsto S_t u\in E$ (called the trajectory of  $u$  
            or/and of $(S_t)_{t\ge 0}$) is continuous.
\end{itemize}
The {\it infinitesimal generator} of the $(C_0)$-semigroup $(S_t)_{t\ge 0}\subset L(E)$ is the linear operator  $A$ from $E$ into $E$  with domain $D(A)$ such that
\begin{align*}
D(A)&= \bigg\{u\in E:\lim_{t\downarrow 0}\frac1t(S_tu-u)\mbox{ exists in the topology of }  E \bigg\},\\
Au&=\lim_{t\downarrow 0}\frac1t(S_tu-u)\ \quad  \mbox{for }u\in D(A).
\end{align*}
Notice that $D(A)$ is dense in~$E$, and $D(A)=E$ if and only if all the trajectories of the  $(C_0)$-semigroup belong to $C^\infty([0,\infty\[;E)$. The name ``generator'' is justified by the fact that if a $(C_0)$-semigroup $(S_t)_{t\ge 0}\subset L(E)$ is locally equicontinuous, then it is uniquely determined by~$A$. (See the proof of the uniqueness theorem in Section~\ref{sec2}.) For one-parameter groups of linear operators things are similar.

\section{The one-parameter group $(G_t)_{t\in\bbR}\subset L((Z'_n)^m)$ generated by $P(\partial/\partial x)$}\label{sec2}

In the present section the matrices $A_{\alpha}\in M_m$, $|\alpha|\le p$, are arbitrary. This is important for the proof of necessity of the Petrovski\u\i\  correctness condition in Theorem 1(iii) of Section~\ref{sec3}.

Since $\phi(t,\xi):=\exp(t\wid P(\xi))$ satisfies the differential equation
$\frac d{dt}\phi=\wid P(\xi)\phi$, the theorem on differentiation of a solution of an ODE  with respect to a parameter
(\cite{H}, Sec.~V.4, Corollary~4.1) implies that $\phi\in C^\infty(\bbR^{n+1};M_m)$. This conclusion may also be (not very easily) obtained by term by term differentiation of the series $\exp(t\wid P(\xi))=\1+\sum_{k=1}^\infty\frac{t^k}{k!}\wid P(\xi)^k$.
Let $\calD'(\bbR^n)$ be the space of distributions on $\bbR^n$
endowed with the topology of uniform convergence on bounded subsets of $C_0^\infty(\bbR^n)$. 
For every $T\in\calD'(\bbR^n)$
the mapping
$C^\infty(\bbR^n)\ni \varphi\mapsto\varphi T\in\calD'(\bbR^n)$
is continuous. Consequently, the formula 
$$
\wid G_tf=e^{t\wid P}f,\ \quad t\in\bbR,\,f\in(\calD'(\bbR^n))^m,
$$
determines a one-parameter group $(\wid G_t)_{t\in\bbR}\subset L((\calD'(\bbR^n))^m)$ of class $(C_0)$ all of whose trajectories belong to $C^\infty(\bbR;(\calD'(\bbR^n))^m)$. 
The infinitesimal generator of this one-parameter group is the multiplication operator 
$\wid P|_{(\calD'(\bbR^n))^m}\in L((\calD'(\bbR^n))^m)$. 
It is easy to
 prove that the one-parameter group 
 $(\wid G_t)_{t\in\bbR}\subset L((\calD'(\bbR^n))^m)$
is locally equicontinuous, i.e. for every compact $K\subset \bbR$ the family of operators $\{\wid G_t:t\in K\}\subset L((\calD'(\bbR^n))^m)$ 
is equicontinuous.

Let  $\calF$  be the $n$-dimensional Fourier transformation defined by
$$
(\calF\varphi)(x)=\int_{\bbR^n}e^{i\<x,\xi>}\varphi(\xi)\,d\xi\ \quad\mbox{for }
\varphi\in \calS(\bbR)\mbox{ and }x\in\bbR^n
$$   
where $\<x,\xi>=\sum_{k=1}^nx_k\xi_k$. Then $\calF$  is an automorphism of $\calS(\bbR^n)$ with inverse $\calF^{-1}=(2\pi)^{-n}\calF^{\lor}$ where $\lor$ denotes the reflection 
in~$0$. Define $Z_n:=\calF C_c^\infty(\bbR^n)$. Then $Z_n=\calF C_c^\infty(\bbR^n)
=\calF^{-1} C_c^\infty(\bbR^n)\subset\calS(\bbR^n)$. 
The space $Z_n$ consists of those functions belonging to $\calS(\bbR^n)$ which have holomorphic extension onto $\bbC^n$ with growth properties characterized by the Paley--Wiener--Schwartz theorem (see \cite{H2}, Theorem 7.3.1). The topology in $Z_n$ is that transported by  $\calF$  from $C_c^\infty(\bbR^n)$. 
The restriction $\calF|_{C_c^\infty(\bbR^n)}$ is a topological isomorphism of $C_c^\infty(\bbR^n)$
onto~$Z_n$, and  $\calF|_{Z_n}=(2\pi)^n\calF^{-1\lor}|_{Z_n}$ is a topological isomorphism of $Z_n$ onto $C_c^\infty(\bbR^n)$. The space  $Z'_n$ is defined as the strong dual 
of~$Z_n$, and the dual mapping of $\calF|_{Z_n}:Z_n\to C_c^\infty(\bbR^n)$ is an isomorphism of $\calD'(\bbR^n)$ onto~$Z'_n$. This last isomorphism {\it extends}
$\calF|_{C_c^\infty(\bbR^n)}$ and for this reason it is still denoted by~$\calF$. It follows that
$$
\<\calF T,u>_{Z'_n\times Z_n}=\<T,\calF u>_{\calD'(\bbR^n)\times C_c^\infty(\bbR^n)}\quad
\mbox{for every }
T\in \calD'(\bbR^n)\mbox{ and }u\in Z_n,
$$
or, what is the same,
$$\displaylines{\indent
\<\calF T,\calF \varphi>_{Z'_n\times Z_n}
=(2\pi)^n\<T,\varphi^{\lor}>_{\calD'(\bbR^n)\times C_c^\infty(\bbR^n)}\hfill\cr
\hfill\quad
\mbox{for every }T\in\calD'(\bbR^n)\mbox{ and }
\varphi\in\calC_c^\infty(\bbR^n).\indent\cr}
$$
For every $k=1,\ldots,n$ the operator $\partial/\partial x_k:Z_n'\to Z_n'$ is defined as dual to the operator  
$-\partial/\partial x:Z_n\to Z_n$, so that $(\partial/\partial x_k)
\calF T=\calF(i\xi_k\cdot T)$ for every 
$T\in\calD'(\bbR^n).$

The actions of  $\calF$  on $(\calD'(\bbR^n))^m$ and of $\calF^{-1}$ on $(Z'_n)^m$
are coordinatewise. From our assertions concerning the one-parameter group 
$(\wid G_t)_{t\in\bbR}$ it follows that 
{\it the formula}
\begin{equation}\label{eq2.1} 
G_t=\calF e^{t\wid P}\calF^{-1},\ \quad t\in\bbR,
\end{equation}
{\it determines a locally equicontinuous one-parameter group
$(G_t)_{t\in\bbR}\subset L((Z'_n)^m)$
of class $(C_0)$ all of whose trajectories belong to $C^\infty(\bbR;(Z'_n)^m)$
and whose infinitesimal generator is} 
$\calF\wid P\calF^{-1}|_{(Z'_n)^m}=P(\partial/\partial x)|_{(Z'_n)^m}$.

If $u_0\in(Z'_n)^m$ and $u(t)=G_tu_0$ for $t\in\bbR$, then $u(\bullet)\in C^\infty(\bbR;(Z'_n)^m)$ and $u(\bullet)$
is a solution of the Cauchy problem \eqref{eq1.4}. The subsequent theorem shows that this Cauchy problem has no other $(Z'_n)^m$-valued solutions.

\medskip
\noindent
\textbf{Uniqueness Theorem.}
{\it Let $t_0\in\]0,\infty]$ and let $I$ be equal to either $[0,t_0\[$ or $\]-t_0,0]$. If the map $I\ni t\mapsto u(t)\in (Z'_n)^m$ belongs to 
$C^1(I;(Z'_n)^m)$ and
$$
\frac{du(t)}{dt}=P\bigg(\frac{\partial}{\partial x}\bigg)u(t)\ \quad
\mbox{for every }t\in I,
$$
then}
$$
u(t)=G_tu(0)\ \quad\mbox{\it for every }t\in I.
$$

\begin{pf}We will prove this theorem for $I=[0,t_0\[$, the proof for $I=\]-t_0,0]$
being similar. Fix any $t\in\]0,t_0\[$ and let $\tau\in[0,t]$. Then
$$
\lim_{h\to0}G_{t-\tau}\frac1h(G_{-h}u(\tau)-u(\tau))
=-G_{t-\tau}P\bigg(\frac{\partial}{\partial x}\bigg)u(\tau)
$$ 
with limit in the topology in $(Z'_n)^m$. Furthermore, by local equicontinuity of the one-parameter group $(G_t)_{t\in\bbR}\subset L((Z'_n)^m)$, the map $\bbR\times(Z'_n)^m
\ni(t,u)\mapsto G_tu\in(Z'_n)^m$ is continuous, and so
$$ 
\lim_{[-\tau,t-\tau]\ni h\to0}G_{t-\tau-h}\frac1h(u(\tau+h)-u(\tau))
=G_{t-\tau}P\bigg(\frac{\partial}{\partial x}\bigg)u(\tau).
$$ 
Consequently,
$$ 
\lim_{[-\tau,t-\tau]\ni h\to0}\frac1h[G_{t-\tau-h}u(\tau+h)-G_{t-\tau}u(\tau)]
=0.
$$
This shows that for every $t\in\]0,t_0]$ the function
$[0,t]\ni\tau\mapsto G_{t-\tau}u(\tau)\in(Z'_n)^m$
has derivative vanishing everywhere on $[0,t]$ (the derivative at the ends of
$[0,t]$ being one-sided). Consequently, $\frac{d}{d\tau}[G_{t-\tau}u(\tau)](\varphi)=0$
for every $\tau\in[0,t]$ and 
$\varphi\in(Z_n)^m$, whence $[u(t)-G_tu(0)](\varphi)
=[G_{t-\tau}u(\tau)](\varphi)|^{\tau=t}_{\tau=0}=0$, and so $u(t)=G_tu(0)$.

Notice that the above argument resembles one used in the proof of E.~R.~van Kampen's uniqueness theorem for solutions of ordinary differential equations. See \cite{K} and \cite{H}, Sec.~III.7.
\end{pf}

\noindent{\bf Remark.}
For every $t\in\bbR$ one has $G_t(Z_n)^m\subset(Z_n)^m$ and $G_t(\calF\calE'(\bbR^n))^m
\subset(\calF\calE'(\bbR^n))^m$. The restricted operators $G_t$ constitute one-parameter $(C_0)$-groups $(G_t|_{(Z_n)^m})_{t\in\bbR}\subset L((Z_n)^m)$
and $(G_t|_{(\calF\calE'(\bbR^n))^m})_{t\in\bbR}\subset L((\calF\calE'(\bbR^n))^m)$
having properties analogous to those of $(G_t)_{t\in\bbR}\subset L((Z'_n)^m)$. Similarly to~$Z_n$, the space $\calF\calE'(\bbR^n)$ has a direct analytical characterization: its elements are those functions which belong to $\calO_M(\bbR^n)$ (the space of slowly increasing $C^\infty$-functions on $\bbR^n$) and have holomorphic extensions onto $\bbC^n$
with growth properties characterized by the Paley--Wiener--Schwartz theorem (\cite{S}, Chapter VII, Theorem XVI; \cite{H2}, Theorem 7.3.1).

\section{The main result}\label{sec3}
 
{\spaceskip.33em plus.22em minus.17em As in Section \ref{sec1}, take a polynomial $P(X)\kern-2.4pt=\kern-2.4pt P(X_1,\ldots,X_n)
\kern-2.4pt\in\kern-2.4pt M_m
[X_1,\ldots,X_n]$.} Let $P(\partial/\partial x)$ be the corresponding matricial partial differential operator with constant coefficients, and let 
$\wid P=\wid P(\xi)= P(i\xi)$
be the symbol of $P(\partial/\partial x)$. Define the {\it stability index} $\omega_0$  of $P(\partial/\partial x)$ by the formula
\begin{equation}\label{eq3.1} 
\omega_0=\sup\{\hRe\lambda:\lambda\in\sigma(\wid P(\xi)),\,\xi\in\bbR^n\}
\end{equation}
where $\sigma(\wid P(\xi))$ denotes the spectrum of the matrix $\wid P(\xi)\in M_m$.

Let  $E$  denote one of the following spaces of $\bbC^m$-valued functions or distributions on $\bbR^n$:
\begin{itemize}{\rm(iii)}
\item[\rm(i)] $E=(\calS(\bbR^n))^m$ where $\calS(\bbR^n)$ is the L.~Schwartz space of rapidly decreasing infinitely differentiable functions on $\bbR^n$,
\item[\rm(ii)] $E=(\calS'(\bbR^n))^m$ where $\calS'(\bbR^n)$ is the L.~Schwartz space of tempered distributions on $\bbR^n$ equipped with the topology of uniform convergence on bounded subsets of $\calS(\bbR^n)$,
\item[\rm(iii)] $E=(C_b^\infty(\bbR^n))^m$ where $C_b^\infty(\bbR^n)$ is the space of those bounded infinitely differentiable functions on $\bbR^n$ whose partial derivatives are all bounded on $\bbR^n$,
\item[\rm(iv)] $E=(H^\infty(\bbR^n))^m$ where $H^\infty(\bbR^n)$ consists of those infinitely differentiable functions on $\bbR^n$ which belong to $L^2(\bbR^n)$ together with all their partial derivatives,
\item[\rm(v)] $E=\{u\in(L^2(\bbR^n))^m:(P(\partial/\partial x))^ku\in(L^2(\bbR^n))^m$
for $k=1,2,\ldots\}$ where partial derivatives are meant in the sense of distributions and the topology of  $E$  is determined by the system of seminorms $\|u\|_k=\|(P(\partial/\partial x))^ku\|_{(L^2(\bbR^n))^m}$, $k=0,1,\ldots.$
\end{itemize}
\medskip

\begin{theor}\label{th1}
Let $P(\partial/\partial x)$ be the matricial partial  differential operator with constant coefficients corresponding to a polynomial $P(X)\in M_m[X_1,\ldots,X_n]$ and let
$\omega_0$ be the stability index of $P(\partial/\partial x)$. Fix whichever of the five spaces  $E$  listed above. Then the following two conditions are equivalent:
\begin{itemize}{\rm(b)}
\item[\rm(a)] $\omega_0<\infty$,
\item[\rm(b)]   the operator $P(\partial/\partial x)|_E$ is the infinitesimal generator of a $(C_0)$-semigroup $(S_t)_{t\ge 0}\subset L(E)$.
\end{itemize}
If $\omega_0<\infty$, then the $(C_0)$-semigroup $(S_t)_{t\ge 0}\subset L(E)$ occurring in {\rm(b)} is unique and $S_t=G_t|_E$
for every $t\ge 0$, where $(G_t)_{t\in \bbR}\subset L((Z_n')^m)$ is the
$(C_0)$-group from Section~{\rm\ref{sec2}}. Furthermore, if $\omega_0<\infty$, then the semigroup $(S_t)_{t\ge 0}$ is of exponential type, i.e.  for every real sufficiently large
$\omega$ the semigroup
$(e^{-\omega t}S_t)_{t\ge 0}\subset L(E)$ is equicontinuous, and
$$
\omega_0=\omega_E
$$
where
$$
\omega_E=\inf\{ \omega\in\bbR:\mbox{the semigroup }(e^{-\omega t}S_t)_{t\ge 0}\subset L(E)\mbox{ is equicontinuous}\}.
$$
\end{theor}

\noindent{\bf Remarks.} 1. In each of the cases (i)--(v) the space $E$ is continuously imbedded in $(Z'_n)^m$ so that if $\omega_0<\infty$, then the equality $S_t=G_t|_E$ and the uniqueness of the
$(C_0)$-semigroup $(S_t)_{t\ge0}\subset L(E)$ generated by $P(\partial/\partial x)|_E$
are consequences of (b) and of the uniqueness theorem from Section~\ref{sec2}. 

2. In case (iii) the equivalence (a)$\Leftrightarrow$(b) follows from results of I.~G.~Petrovski\u\i\  \cite{P} not involving one-parameter semigroups. Condition (a), now called the {\it Petrovski\u\i\  correctness condition}, was used in \cite{P} in a seemingly weaker form which was later  proved to be equivalent to (a), according to a conjecture formulated in~\cite{P}. From Theorem 1(i) it follows that if (a) holds, then $S_t=\calF e^{t\wid P}\calF^{-1}=T_t\,*$ for every $t\in[0,\infty\[$ where $T_t\in(\calS'(\bbR^n))^m$. Formulas (8.5)--(8.6) in Sec.~7.8 of A.~Friedman's book \cite{F} exhibit the structure of distributions
$T_t$ and yield some results of type (a)$\Rightarrow$(b) not related to exponential semigroups, generalizing the above-mentioned results of
Petrovski\u\i.
 
3. One may call $\omega_E$ the {\it characteristic exponent} or the 
{\it equicontinuity index} of the semigroup $(S_t)_{t\ge0}\subset L(E)$. A semigroup $(S_t)_{t\ge0}\subset L(E)$ with $\omega_E$ finite may be called an {\it exponential semigroup}. The exponential
$(C_0)$-semigroups in an l.c.v.s.~reduce (by multiplication by real exponential function of the parameter $t$) to equicontinuous 
$(C_0)$-semigroups for which a theory of Hille--Yosida type  is presented in Chapter IX of the monograph of K.~Yosida~\cite{Y}. In a Banach space all one-parameter $(C_0)$-semigroups are exponential, but a similar statement is not true for non-normed spaces. See for instance \cite{Ki},~p.~4.

4. The idea of using one-parameter semigroups in connection with Cauchy's problem \eqref{eq1.4} is taken from papers of G.~Birkhoff, T.~Mullikin and T.~Ushijima \cite{B-M}, \cite{B} and \cite{U}, and from Section I.8 of S.~G.~Krein's monograph \cite{Kr}. Let us stress that in \cite{B} the equality $\omega_0=\omega_E$ is discussed.

5.  A  result similar to Theorem 1(i) concerning Cauchy's problem for  the equation
$$
P\bigg(\frac{\partial}{\partial t},\frac{\partial}{\partial x_1},\ldots,
\frac{\partial}{\partial x_n}\bigg)=0
$$
where $P(X_0,X_1,\ldots,X_n)=\sum_{k=0}^mP_k(X_1,\ldots,X_n)X_0^k$
is a scalar polynomial of $n+1$ variables is stated in the book of J.~Rauch \cite{R} as Theorem 2 on p.~128. Since the polynomial
$P_m(X_1,\ldots,X_n)$ need not reduce to a constant, the theorem of Rauch does not follow from Theorem 1(i) (or vice versa).

6. The spaces  $E$  in cases (i)--(iv) are standard, not depending on $P(\partial/\partial x)$. The space  $E$  in case~(v), depending on
$P(\partial/\partial x)$, was introduced by T.~Ushijima \cite{U} who proved the equivalence (a)$\Leftrightarrow$(b) in this case. Let  
$X\kern-2.5pt=\kern-2.5pt(L^2(\bbR^n))^m$, and let  $A$  be the operator from  $X$ into $X$  with domain $D(A)$ such that $D(A)=\{u\in X:P(\partial/\partial x)
u\in X\}$
and $Au=P(\partial/\partial x)u$  for $u\in D(A)$. Then for the space $E$  of Theorem 1(v) one has $E=D(A^\infty):=\bigcap_{k=1}^\infty D(A^k)$ and, in the terminology of \cite{U}, assertion (b) of Theorem 1(v) means that the operator  $A$  is $D(A^\infty)$-{\it well posed}. The $D(A^\infty)$-well posedness of an operator  $A$  from a Banach space into itself is one of the central notions of the ACP-theoretical paper~\cite{U}. In Sec.~1.3 of \cite{Ki} the notion of $D(A^\infty)$-well posedness  is also elucidated by some facts not mentioned in~\cite{U}.

7. Results similar to ${\rm(a)\Rightarrow(b)}\wedge(\omega_E\le\omega_0)$ of Theorems 1(iii) and 1(iv) constitute a part of Theorem 4.1 of \cite{H-H-N} obtained by means of one-parameter regularized semigroups of operators and Fourier multipliers.

\section{Banach and Hilbert spaces adapted to $P(\partial/\partial x) $. Constructions of 
G.~Birkhoff and  S.~D.~Eidelman--S.~G.~Krein}\label{sec4}           

In Theorem 1 all the spaces  $E$  are not normed. In order to compare Theorem 1 with earlier results related to some special Banach and Hilbert spaces  $E$  depending on $P(\partial/\partial x) $ let us recall two constructions:
\begin{itemize}{\rm(vii)}
\item[\rm(vi)] $E=\calB_{\calN,p}$ {where $p\in[1,\infty\[$ and
$\calB_{\calN,p}$ is the Banach space introduced by G.~Birkhoff and T.~Mullikin in {\rm\cite{B-M}} and~{\rm\cite{B}},}
\item[\rm(vii)]$E=\calL_B$  {where $\calL_B$ is the Hilbert space constructed by S.~D.~Eidelman and
 S.~G.~Krein.}
\end{itemize}
$\calB_{\calN,p}$ depends upon the K.~Baker map $\calN$ whose existence was proved in \cite{Ba}, and which is a Borel measurable map 
$\calN:\bbR^n\to M_m$  such that for every $\xi\in\bbR^n$ the matrix  $\calN(\xi)$ is invertible and $\calN(\xi)\wid P(\xi)\calN(\xi)^{-1}$   is a Jordan matrix whose diagonal elements belong to $\sigma(\wid P(\xi))$, directly over-diagonal elements are equal to zero or one, and all other elements are equal to zero. Given a Baker map
$\calN$ one defines $\calB_{\calN,p}$ as a linear subset of
$(Z'_n)^m$ consisting of all those elements  $u$  of
$(Z'_n)^m$ for which the distribution $\calF^{-1}u\in(\calD'(\bbR^n))^m$
is a $\bbC^m$-valued Lebesgue measurable function on $\bbR^n$ such that
\begin{equation}\label{eq4.1}
\|u\|_{\calN,p}=\bigg(\int_{\bbR^n}\|\calN(\xi)(\calF^{-1}u)(\xi)\|_{\bbC^m}^p\,d\xi\bigg)^{1/p}<\infty.
\end{equation}
Equipped with the norm $\|\|_{\calN,p},\calB_{\calN,p}$ is a Banach space continuously imbedded in $(Z'_n)^m$, and $\calB_{\calN,2}$ is a Hilbert space. If the Baker function $\calN:\bbR^n\to M_m$ is bounded on
$\bbR^n$ and $p\in[2,\infty\[$, then, by the Hausdorff--Young theorem, $\calB_{\calN,p}\supset (L^q(\bbR^n))^m$, where $q=p/(p-1)\in\]1,2]$.
 
Construction of $\calL_B$, described in Section I.8 of the monograph \cite{Kr}, requires the assumption that $\omega_0<\infty$. If $\omega_0<\infty$, then $\calL_B$ is the domain of the selfadjoint strictly positive definite square root of some selfadjoint strictly positive definite operator $B(\partial/\partial x)$ acting in $(L^2(\bbR^n))^m$. For any given $\omega_1\in\]\omega_0,\infty\[$ the operator $B(\partial/\partial x)$
may be constructed as a matricial partial differential operator with constant coefficients whose symbol $\wid B$ has the following property:
\medskip

{\it for every $\xi\in\bbR^n$ the matrix $\wid B(\xi)\in M_m$ is hermitian such that}
\begin{equation}\label{eq4.2}
\wid B(\xi)\ge\1\quad\mbox{\it and}\quad \wid B(\xi)\wid P(\xi)
+\wid P(\xi)^*\wid B(\xi)\le 2\omega_1\wid B(\xi).
\end{equation}
Let $\calN(\xi)$ be the hermitian strictly positive definite square root of $\wid B(\xi)$. Then $\calN(\xi)=\frac1{2\pi i}\int_C z^{1/2}
(z\1-\wid B(\xi))^{-1}\,dz$
where  $C$ is a closed rectifiable curve contained in
$\{z\in \bbC:\hRe z\ge 1/2\}$ and winding once about $\sigma(\wid B(\xi))$, which is a finite subset of $[1,\infty\[$. 
It follows that $\calN(\xi)$ is a $C^\infty$-function of~$\xi$. By \eqref{eq4.2} for every $\xi\in\bbR^n$ one has
$$
\|\calN(\xi)^{-1}\|_{L(\bbC^m)}\le 1\quad\mbox{and}\quad\calN(\xi)
\wid P(\xi)\calN(\xi)^{-1}+(\calN(\xi)\wid P(\xi)\calN(\xi)^{-1})^*
\le2 \omega_1\1.
$$
The norm in ${\calL}_B$ is defined by the formula
\begin{align}\label{eq4.3}   
\|u \|_{{\calL}_B}
&=\bigg\|\bigg(B\bigg(\frac{\partial}{\partial x}\bigg)\bigg)^{1/2}u \bigg\|_{(L^2(\bbR^n))^m}\notag\\
&=(2\pi)^n\bigg(\int_{\bbR^n}\|\calN(\xi)(\calF^{-1}u)(\xi)\|^2_{\bbC^m}
\,d\xi\bigg)^{1/2}. 
\end{align}
${\calL}_B$ is a Hilbert space continuously imbedded in $(L^2(\bbR^n))^m$. 
\goodbreak

\begin{theor}\ 
\begin{itemize}{\rm(II)}
\item[{\rm(I)}]  Suppose that 
$p\in[1,\infty\[$, $\calN$ is a Baker function for $\widetilde P$, and 
$\calN\in\break L_{\rm loc}^p(\bbR^n;M_m)$. 
Then the same statements as in Theorem $1$ are true for $E=\calB_{\calN,p}$ provided
$P(\partial/\partial x)|_E$ is replaced by the restriction of $P(\partial/\partial x)$ to the set $\{u\in \calB_{\calN,p}:P(\partial/\partial x)u\in\calB_{\calN,p}\}$.
\item[\rm(II)]   If $\omega_0<\infty$, $\omega_1\in\]\omega_0,\infty\[$, $E={\calL}_B$ is constructed so that {\rm\eqref{eq4.2}} is satisfied, and $P(\partial/\partial x)|_E$
is replaced by the restriction of $P(\partial/\partial x)$ to the set
${\{u\in {\calL}_B:P(\partial/\partial x)u\in {\calL}_B\}}$, then {\rm(b)} holds and $\omega_0\le\omega_E\le\omega_1$.
\end{itemize}
\end{theor}

The implication ${\rm(a)\Rightarrow(b)}\wedge(\omega_0=\omega_{\calB_{\calN,p}})$ of item (I) was proved by G.~Birk\-hoff 
\cite{B}  without the assumption that the Baker function is locally bounded or locally integrable. The statement (II), except the inequality $\omega_0\le\omega_{{\calL}_B}$, goes back to S.~D.~Eidelman and S.~G.~Krein. The proof is given in Section I.8 of S.~G.~Krein's monograph \cite{Kr}. Apart from the results proved in \cite{B} and \cite{Kr}, Theorem 2 contains information that in both cases,
$E={\calL}_B$ and $E=\calB_{\calN,p}$ with some special  $\calN$, condition (b) implies that $\omega_0\le\omega_E$. 

In order to prove this last statement  suppose that (b) holds, and let $(S_t)_{t\ge0}\subset L(E)$  be the (unique) semigroup of class $(C_0)$
occurring in (b). The equicontinuity index $\omega_E$  of this semigroup is finite and equal to the characteristic exponent
$\lim_{t\to\infty}\frac1t\log\|S_t\|_{L(E)}$ of the function $t\mapsto\|S_t\|_{L(E)}$. Consequently, by \eqref{eq4.1} and \eqref{eq4.3}, for every $\varepsilon>0$ there is $K_{\varepsilon}\in[1,\infty\[$ such that
$$
\bigg(\!\int_{\bbR^n}\kern-2.25pt\|\calN(\eta)(\calF^{-1}S_t\calF\varphi)
(\eta)\|_{\bbC^m}^p\,d\eta \bigg)^{1/p}
\!\!\le\!
K_{\varepsilon}e^{(\omega_E+\varepsilon)t}
\bigg(\!\int_{\bbR^n}\kern-2.25pt\|\calN(\eta)\varphi(\eta)\|_{\bbC^m}^p\,d\eta \bigg)
^{1/p}
$$
for every $t\in[0,\infty\[$ and $\varphi\in(L_{c}^\infty(\bbR^n))^m$, i.e. $M_m$-valued $\varphi$ which is bounded and measurable on
$\bbR^n$ and has compact support. By the uniqueness theorem from Section~\ref{sec2}, one has $S_t=G_t|_E$, so that $(\calF^{-1}S_t\calF\varphi)(\eta)=e^{t\wid P(\eta)}\varphi(\eta)$ and hence
\begin{multline}\label{eq4.4}
\bigg(\int_{\bbR^n}\|\calN(\eta)e^{t\wid P(\eta)}\varphi(\eta)\|_{\bbC^m}^p
\,d\eta\bigg)^{1/p}\\
\le
K_{\varepsilon}e^{(\omega_E+\varepsilon)t}
\bigg(\int_{\bbR^n}\|\calN(\eta)\varphi(\eta)\|_{\bbC^m}^p\,d\eta\bigg)^{1/p}
\end{multline}
whenever $t\in[0,\infty\[$ and $\varphi\in(L_c^\infty(\bbR^n))^m$. 

For any $t\in[0,\infty\[$ and $\xi\in\bbR^n$ pick a $z_{t,\xi}\in
\bbC^m$ such that $\|z_{t,\xi}\|_{\bbC^m}=1$ and 
\begin{equation}\label{eq4.5}
\|\calN(\xi)e^{t\wid P(\xi)}\|_{L(\bbC^m)}
=\|\calN(\xi)e^{t\wid P(\xi)}z_{t,\xi}\|_{\bbC^m}.
\end{equation}
For every $\xi\in\bbR^n$ and $ r>0$ let $\phi_{\xi,r}\in L_{c}^\infty(\bbR^n)$  be a non-negative function such that $\int_{\bbR^n}\phi_{\xi,r}
(\eta)\,d\eta=1$, the support of $\phi_{\xi,r}$ is equal to the ball with center at  $\xi$  and radius  $r$, and $\phi_{\xi,r}$ is constant in this ball. 
Since, in case (I), $\calN\in L_{\rm loc}^p(\bbR^n;M_m)$, it follows that there
is a set $\Z_t\subset\bbR^n$ of $n$-dimensional Lebesgue measure zero such that if $\xi\in\bbR^n\setminus \Z_t$, then $\xi$ is a Lebesgue point of the 
locally integrable function $\|\calN(\cdot)e^{t\widetilde P(\cdot)}\|^p$, i.e.
\begin{equation}\label{eq4.6}
\lim_{r\downarrow 0}\int_{\bbR^n}
|\|\calN(\eta)e^{t\wid P(\eta)}\|^p
-\|\calN(\xi)e^{t\wid P(\xi)}\|^p|
\phi_{\xi,r}(\eta)\,d\eta=0.
\end{equation}
See \cite{S-K}, Theorem 5.3, p.~164. 
In case (II) this difficult theorem need not be used because (4.6) with
$\Z_t=\emptyset$ holds by virtue of continuity of~$\calN$.
From \eqref{eq4.5} and \eqref{eq4.6} it follows that whenever 
$t\in[0,\infty\[$ and $\xi\in\bbR^n\setminus \Z_t$, 
 then
\begin{equation*}
\|\calN(\xi)e^{t\wid P(\xi)}\|
=\lim_{r\downarrow 0}\bigg(\int_{\bbR^n}\|\calN(\eta)e^{t\wid P(\eta)}
(\phi_{\eta,r}(\eta))^{1/p}z_{t,\xi}\|^p\,d\eta\bigg)^{1/p}.
\end{equation*} 
Hence, applying \eqref{eq4.4} to $\varphi(\eta)=(\phi_{\xi,r}(\eta))^{1/p}z_{t,\xi}$, one concludes that whenever  $t\in[0,\infty\[$
 and $\xi\in\bbR^n\setminus (\Z_0\cup\Z_t)$, then
$$\displaylines{
\|\calN(\xi)e^{t\wid P(\xi)}\|\hfill\cr\hfill
\le K_{\varepsilon} e^{(\omega_E+\varepsilon)t}
\lim_{r\downarrow 0}\bigg(\int_{\bbR^n}\|\calN(\eta)(\phi_{\xi,r}(\eta))^{1/p}
z_{t,\xi}\|^p\,d\eta\bigg)^{1/p}
\le K_{\varepsilon}e^{(\omega_E+\varepsilon)t}
\|\calN(\xi)\|\cr}
$$
and consequently
$$
\|e^{t\wid P(\xi)} \|
\le K_{\varepsilon}\|\calN(\xi)\|\,\|\calN(\xi)^{-1}
\|e^{(\omega_E+\varepsilon)t}. 
$$
By Proposition 2.2, p.~251, and Corollary 2.4, p.~252, in \cite{E-N} it follows that
$$
\max\{\hRe\lambda:\lambda\in\sigma(\wid P(\xi))\}
=\lim_{t\to\infty}\frac1t\log\|e^{t\wid P(\xi)}\|
\le\omega_E+\varepsilon
$$
for every $\xi\in\bbR^n\setminus \bigcup_{k\in\bbN_0}\Z_k$.  Since  $\bigcup_{k\in\bbN_0}\Z_k$  has measure zero and
$\max\{\hRe\lambda:\lambda\in\sigma(\wid P(\xi))\}$ depends continuously on $\xi$, one concludes that $\omega_0\le\omega_E+\varepsilon$. This implies the inequality
$\omega_0\le\omega_E$ because $\varepsilon>0$ is arbitrary.

\section{G\aa rding's lemma}\label{sec5}
 
\medskip
\noindent{\bf Theorem} {\rm(L.~G\aa rding)}. {\it Consider a polynomial of  $n+1$
variables with complex coefficients
$$
p(X_0,X_1,\ldots,X_n)\!=\!\sum_{k=1}^{m}p_k(X_1,\ldots,X_n)X_0^k+p_0(X_1,\ldots,X_n)
\in\bbC[X_0,\ldots,X_n].
$$
Suppose that $p_m(\xi_1,\ldots,\xi_n)>0$ for every $(\xi_1,\ldots,\xi_n)\in\bbR^n$. For every $r\in[0,\infty\[$ define
\begin{multline*}\indent
\Lambda(r)=\sup\{\hRe\lambda:\lambda\in\bbC,\, (\xi_1,\ldots,\xi_n)
\in\bbR^n,\\
\,\|(\xi_1,\ldots,\xi_n)\|\le r,\,p(\lambda,\xi_1,\ldots,\xi_n)= 0\}.\indent
\end{multline*}
Then there is a real  $A$  and a rational $\alpha$ such that}
$$
\Lambda(r)=Ar^\alpha(1+o(1))\ \quad\mbox{\it as }r\to\infty.
$$

This theorem was formulated by L.~G\aa rding in \cite{G} as the Lemma on p.~11. The argument consisted in
\begin{itemize}{\rm(B)} 
\item[\rm(A)] explaining the existence of a polynomial $q(Y_1,Y_2)\in\bbC[Y_1,Y_2]$ such that $q(r,\Lambda(r))=0$
for every  $r\in[0,\infty\[$,
\item[\rm(B)]  applying the Puiseux series expansions of algebraic functions $\calR$ of one 	complex variable  $z$  satisfying the equation $q(z,\calR(z))=0$.
\end{itemize}

L.~H\"ormander \cite{H1}, proof of Lemma 3.9, noticed that stage (A) may be realized by an application of a theorem asserting that the projection onto $\bbR^d$ of a semi-algebraic subset of $\bbR^{d+n}$  is a semi-algebraic subset of $\bbR^d$. This projection theorem
may be proved by an argument similar to that from
A.~Seidenberg's proof \cite{Se} of the decision theorem  of A.~Tarski
(belonging to mathematical logic). Detailed presentations of 
Seidenberg's proof in the case of the projection theorem are given in \cite{G2} and \cite{F}. 
L.~H\"ormander's proof of the
 projection theorem presented in the Appendix to \cite{H2} is based on an argument resembling that from P.~J.~Cohen's proof [C1, 2] of the decision theorem. 
In [Ki] the reasonings of stage (B) of G\aa rding's proof are presented with exact references to the theory of algebraic functions of one complex variable presented in \cite{S-Z}.

The above theorem of G\aa rding yields at once

\begin{coro}  Let $P(X_1,\ldots,X_n)\in M_m[X_1,\ldots,X_n]$ and suppose that there is
$C\in\]0,\infty\[$ such that
$$
\max\{\hRe\lambda:\lambda\in\sigma(\wid P(\xi))\}\le C+C\log(1+|\xi|)\ \quad
\mbox{for every }\xi\in\bbR^n.
$$
Then
$$
\sup\{\hRe\lambda:\lambda\in\sigma(\wid P(\xi)),\,\xi\in\bbR^n\}<\infty.
$$
\end{coro}

This corollary was formulated as a conjecture by I.~G.~Petrovski\u\i\  in a footnote on p.~24 of \cite{P}.

\section{Interpolation polynomials and estimations of  $e^{t\wid P(\xi)}$}\label{sec6} 

\medskip

\noindent{\bf Assumption (A).}
 Let $\lambda_1,\ldots,\lambda_m\in\bbC$. Denote by  $S$  the set $\{\lambda_1,\ldots,\lambda_m\}$. For every
$\lambda\in S$ denote by $m(\lambda)$ the number of occurrences of $\lambda$ in the sequence $\lambda_1,\ldots,\lambda_m$. Let 
$p(\lambda)=p_0+p_1\lambda+\cdots+ p_d\lambda^d$
be a polynomial of degree $d$ with complex coefficients. Let  $f$  be a function holomorphic in an open set $O\subset\bbC$
containing~$S$. Choose $r>0$ such that $K:=\bigcup_{\lambda\in S}
\{z\in\bbC:|z-\lambda|\le r\}
\subset O$ and let $\calC $   be the boundary of  $K$  oriented so that
${\rm Index}(\calC,\lambda)=1$ for every ${\lambda\in S}$.

\medskip

\noindent{\bf Theorem I.} {\it  Under assumption {\rm(A)} the following two conditions are equivalent:

\begin{itemize}{\rm(ii)}
\item[\rm(i)] $p^{(k)}(\lambda)=f^{(k)}(\lambda)$ for every $\lambda\in S$ and $k=0,\ldots,m(\lambda)-1$,
\item[\rm(ii)] $p(A):=p_0\1+\sum_{k=1}^dp_kA^k=\frac1{2\pi i}\int_{\calC }f(z)
(z\1-A)^{-1}\,dz$
for every matrix ${A\in M_m}$ such that $\sigma(A)=S$ and for every $\lambda\in\sigma(A)$ the spectral multiplicity of $\lambda$ is equal to~$m(\lambda)$.
\end{itemize}
\looseness2
There is exactly one polynomial  p  of order no greater than
$m-1$ satisfying~{\rm(i)}.}

\medskip

For a polynomial  $p$  of arbitrary degree  $d$  the equivalence (i)$\Leftrightarrow$(ii) may be deduced either from Theorems 5, 8 and 10 of Section VII.1 of \cite{D-S}, or from Theorems 134, 138 and 234 of \cite{G-L}. For $d\le m-1$ the equivalence (i)$\Leftrightarrow$(ii) is a part of Fact 1 stated in \cite{Hig}. In connection with (i),
 interpolation terminology is used:  $p$  is called the {\it interpolation polynomial} for  $f$, the numbers $\lambda\in S$ are called the {\it nodes of interpolation}, and $m(\lambda)$ is the {\it multiplicity of the node}~$\lambda$.

\medskip

\noindent
{\bf Theorem II.}  {\it For a polynomial  p  of degree no greater that $m-1$ conditions {\rm(i)} and {\rm(ii)} are equivalent to either of the following conditions:
\begin{itemize}{\rm(iii)}
\item[\rm(iii)] $p(\lambda)=c_0+c_1(\lambda-\lambda_1)+c_2(\lambda-\lambda_1)(\lambda-\lambda_2)+\cdots$
$$\displaylines{\hfill
{}+c_{m-1}(\lambda-\lambda_1)\cdots(\lambda-\lambda_{m-1})\ \quad
\refstepcounter{equation}(\theequation)\label{eq6.1}
\cr}
$$
where\vg
\begin{align}\label{eq6.2}
c_k&=c_k(f;\lambda_1,\ldots,\lambda_{k+1})\notag\\
&=\frac1{2\pi i}\int_{\calC }f(z)(z-\lambda_1)^{-1}
(z-\lambda_2)^{-1}\cdots(z-\lambda_{k+1})^{-1}dz
\end{align}
for $k=0,\ldots,m-1$,\goodbreak
\item[\rm(iv)] $p(\lambda)\!=\!\disp\frac1{2\pi i}\int_{\calC }f(z)\sum_{\mu=1}^m\frac1{\mu!}
\frac{Q^{(\mu)}(z)}{Q(z)}(\lambda-z)^{\mu-1}\,dz$ where
$Q(z)\!=\!\disp\prod_{k=1}^m({z-\lambda_k})$.
\end{itemize}
If $O\supset{\rm conv} S$, then 
\begin{align}\label{eq6.3}
c_k(f;\lambda_1,\ldots,\lambda_{k+1})
={}&\int_0^1\int_0^{t_1}\cdots\int_0^{t_{k-1}}f^{(k)}((1-t_1)\lambda_1+(t_1-t_2)\lambda_2+\cdots\notag\\
&{}+(t_{k-1}-t_k)\lambda_k+t_k\lambda_{k+1})\,dt_k\,dt_{k-1}\cdots dt_1
\end{align}
for $k=0,\ldots,m-1$.  Furthermore,
\begin{equation}\label{eq6.4}
\frac1{\mu!}\frac{Q^{(\mu)}(z)}{Q(z)}=\tau_\mu\bigg(\frac1{z-\lambda_1},
\ldots,\frac1{z-\lambda_m}\bigg) \ \quad           
\mbox{for }\mu=1,\ldots,m
\end{equation}
where $\tau_\mu(x_1,\ldots,x_m)=\sum_{1\le i_1<\cdots<i_\mu\le m}
x_{i_1}\cdots x_{i_\mu},$ $\mu=1,\ldots,m$, are the elementary symmetric polynomials of  m variables $x_1,\ldots,x_m$. Consequently, condition {\rm(iv)} may be written in the equivalent form:
$$\displaylines{{\rm(iv)'}\
p(\lambda)=a_0+a_1\lambda+\cdots+a_{m-1}\lambda^{m-1}\hfill
\refstepcounter{equation}(\theequation)\label{eq6.5}\cr
\phantom{{\rm(iv)'}}\,
\mbox{where}\hfill\cr
\indent
a_k=a_k(f;\lambda_1,\ldots,\lambda_m)=\sum_{l=0}^{m-k-1}\binom{k+l}{k}I_{k+l+1}^{l}
\ \quad\mbox{for }k=0,\ldots,m-1\hskip7pt\hfill
\refstepcounter{equation}(\theequation)\label{eq6.6}
\cr\noalign{\goodbreak}
\phantom{{\rm(iv)'}}\, 
\mbox{and}\hfill\cr
\phantom{{\rm(iv)'}}\, 
I_{\mu}^{l}=\frac1{2\pi i}\int_{\calC }f(z)(-z)^l\tau_{\mu}
\bigg(\frac1{z-\lambda_1},\ldots,
\frac1{z-\lambda_m}\bigg)\,dz\hfill
\refstepcounter{equation}(\theequation)\label{eq6.7}
\cr
\phantom{{\rm(iv)'}}\, 
\mbox{for $\mu=1$ and  $l=0,\ldots,\mu-1$.}\hfill 
\cr}
$$\vskip-\lastskip}
\vskip12pt
 
The explicit formulas \eqref{eq6.2} and \eqref{eq6.3} for the coefficients $c_k$ in the Newton form \eqref{eq6.1} of the interpolation polynomial of degree no greater than $m-1$ are deduced in Section I.4.2 and I.4.3 of A.~O.~Gelfond's book \cite{Ge}. E.~A.~Gorin \vg \cite{G1} inferred from \eqref{eq6.2} that the coefficients $a_0,\ldots,a_{m-1}$ of the interpolation polynomial in the form \eqref{eq6.5} are linear combinations of the integrals
\begin{equation}\label{eq7.8}
I_{i_1,\ldots,i_k}^l=\frac1{2\pi i}\int_{\calC }f(z)(-z)^l
(z-\lambda_{i_1})^{-1}\cdots(z-\lambda_{i_k})^{-1}\,dz.
\end{equation}
The exact computation of these linear combinations by Gorin's method is possible but troublesome. We will prove \eqref{eq6.6}--\eqref{eq6.7} by another method, based on Theorem~I. Notice that, in connection with Cauchy's problem for a system of PDE with constant coefficients, formulas similar to \eqref{eq6.6}--\eqref{eq6.7} were used by 
E.~A.~Gorin \cite{G1} and T.~Ushijima \cite{U}.

\medskip

\noindent{\bf Proof of (iv) and (\ref{eq6.4})--(\ref{eq6.7}).} Take a matrix $A\in M_m$ such that $\sigma(A)=S$ and for every
$\lambda\in\sigma(A)$ the spectral multiplicity of $\lambda$ is equal to $m(\lambda)$. By Taylor's formula and the Cayley--Hamilton theorem,
$$
Q(z)\1+\sum_{\mu=1}^m\frac1{\mu!}Q^{(\mu)}(z)(A-z\1)^\mu=Q(A)=0,
$$
whence
$$
(z\1-A)^{-1}=\sum_{\mu=1}^m\frac1{\mu!}\frac{Q^{(\mu)}(z)}{Q(z)}
(A-z\1)^{\mu-1}\ \quad\mbox{for every }z\in\bbC\setminus S,
$$
and so the polynomial $p(\lambda)$ occurring in (iv) satisfies (ii). Moreover, its degree  is no greater than $m-1$, so that, by Theorem~I, it is the unique polynomial satisfying (ii). This proves the equivalence (ii$\Leftrightarrow$(iv) in the class of polynomials of degree no greater than $m-1$.

In order to prove \eqref{eq6.4} notice that
$Q(z)=\tau_m(z-\lambda_1,\ldots,z-\lambda_m)$ and
$$
\frac{d}{dz}\tau_\mu(z-\lambda_1,\ldots,z-\lambda_m)
=(m-\mu+1)\tau_{\mu-1}(z-\lambda_1,\ldots,z-\lambda_m)
$$
for $z\in\bbC$ and $\mu=1,\ldots,m$ where $\tau_0\equiv 1$. Consequently,
$$
Q^{(\mu)}(z)=\bigg(\frac{d}{dz}\bigg)^\mu\tau_m(z-\lambda_1,\ldots,z-\lambda_m)=\mu!\tau_{m-\mu}(z-\lambda_1,\ldots,z-\lambda_m),
$$
and so
$$
\frac1{\mu!}\frac{Q^{(\mu)}(z)}{Q(z)}=\frac{\tau_{m-\mu}(z-\lambda_1,\ldots,z-\lambda_m)}{\tau_{m}(z-\lambda_1,\ldots,z-\lambda_m)}
=\tau_\mu\bigg(\frac1{z-\lambda_1},\ldots,\frac1{z-\lambda_m}\bigg)
$$
for $z\in\bbC\setminus S$ and $\mu=1,\ldots,m$. Therefore the polynomial occurring in (iv) may be written in the form 
$$
p(\lambda)=\frac1{2\pi i}\int_{\calC }f(z)\bigg[\sum_{\mu=1}^m\tau_\mu
\bigg(\frac1{z-\lambda_1},\ldots,\frac1{z-\lambda_m}\bigg)
(\lambda-z)^{\mu-1}\bigg]\,dz,
$$
whence the formulas \eqref{eq6.6}--\eqref{eq6.7} for the coefficients
$a_k$, $k=0,1,\ldots,\allowbreak m-\nobreak1$, occurring in \eqref{eq6.5} follow by applying the binomial formula to $(\lambda-\nobreak z)^{\mu-1}$.

I.~M.~Gelfand and G.~E.~Shilov in Sec.~II.6.1 of their book \cite{G-S3}
 have reproduced the proof  of A.~O.~Gelfond's formula \eqref{eq6.3} and observed that this formula implies at once the important inequality
\begin{equation}\label{eq6.9}
\begin{aligned}
\|e^{tA}\|_{L(\bbC^m)}&\le e^{\omega t}\bigg(1+\sum_{k=1}^{m-1}
\frac{(2t)^k}{k!}\|A\|_{L(\bbC^m)}^k\bigg),\\
\omega&=\max\{\hRe\lambda:\lambda\in\sigma(A)\},
\end{aligned}
\end{equation}    
for every $A\in M_m$ and $t\in[0,\infty\[$. The proof of\eqref{eq6.3} and \eqref{eq6.9} is also presented in Sec.~7.2 of A.~Friedman's  book \cite{F}.

The inequality \eqref{eq6.9} is crucial for the proofs of our Theorem 1 from Section~\ref{sec3} in cases (i)--(iv). In case (v) we follow T.~Ushijima \cite{U} and instead of using \eqref{eq6.9} we base on estimation of some Gorin's integrals. This method yields the following

\medskip

\noindent{\bf Proposition.} {\it  Let $\wid P\in M_m[\xi_1,\ldots,\xi_m]$ be the symbol of a matricial differential operator
$P(\partial/\partial x)$ with constant coefficients described in Section~{\rm\ref{sec1}}. Suppose that the Petrovski\u\i\  correctness condition is satisfied:
\begin{equation}\label{eq6.10}
\sup\{\hRe\lambda:\lambda\in\sigma(\wid P(\xi)),\,
\xi\in\bbR^n\}=\omega_0<\infty.
\end{equation}
Then there are functions $p_k\in C^\infty(\bbR^{1+n};\bbC)$, $k=0,1,\ldots,2m$, such that
\begin{equation}\label{eq6.11}
e^{t\wid P(\xi)}=p_0(t,\xi)\1+\sum_{k=1}^{2m}p_k(t,\xi)\wid P(\xi)^k
\ \quad\mbox{for every }(t,\xi)\in\bbR^{1+n}
\end{equation}
and
\begin{equation}\label{eq6.12}
\sup\{e^{-(\omega_0+\varepsilon)t}|p_k(t,\xi)|:t\in[0,\infty\[,\,
\xi\in\bbR^n\}<\infty
\end{equation}
for every $\varepsilon>0$ and $k=0,1,\ldots,2m$.
}                                 

\begin{pf}By Theorems 5 and 10 in Sec.~VII.1 of \cite{D-S}, or by Facts 1 and 8 in Sec.~1 of \cite{Hig}, for every $(t,\xi)\in\bbR^{1+n}$ one has 
\begin{align*}
e^{t\wid P(\xi)}
={}&\lim_{N\to\infty}\bigg(\1+t\wid P(\xi)+\cdots+\frac{t^N}{N!}\wid P(\xi)^N\bigg)\\
={}&\lim_{N\to\infty}\frac1{2\pi i}\int_{\calC_{\xi} }\bigg(1+tz+\cdots+\frac{t^N}{N!}z^n\bigg)(z\1-\wid P(\xi))^{-1}\,dz\\\noalign{\goodbreak}
={}&\frac1{2\pi i}\int_{\calC_\xi }e^{tz}(z\1-\wid P(\xi))^{-1}\,dz\\
={}&\frac1{2\pi i}\int_{\calC_{\xi} }(z-z_0)^{m+1}(z\1-\wid P(\xi))^{-1}\,dz\\
&{}\times\frac1{2\pi i}\int_{\calC_{\xi} }(z-z_0)^{-m-1}e^{tz}(z\1-\wid P(\xi))^{-1}\,dz\\
={}&(\wid P(\xi)-z_0\1)^{m+1}\cdot\frac1{2\pi i}\int_{\calC _{\xi}}(z-z_0)^{-m-1}e^{tz}(z\1-\wid P(\xi))^{-1}\,dz
\end{align*}   
where $z_0\in\bbC$ is a point such that $\hRe z_0>\omega_0$ and
$\calC _{\xi}$ is a rectifiable closed path contained in
$\{z\in \bbC\setminus\sigma(\wid P(\xi)):\hRe z<\hRe z_0\}$ and winding once about $\sigma(\wid P(\xi))$. By Theorem I and Theorem II(iv) and (iv)$'$, it follows that
\begin{align}
e^{t\wid P(\xi)}={}&(\wid P(\xi)-z_0\1)^{m+1}(a_0(t,\xi)\1\notag\\
\label{eq6.13}
&{}+a_1(t,\xi)
\wid P(\xi)+\cdots+a_{m-1}(t,\xi)\wid P(\xi)^{m-1})
\end{align}
for every $(t,\xi)\in\bbR^{1+n}$ where
\begin{align}
a_k(t,\xi)
={}&\sum_{l=0}^{m-1-k}\binom{k+l}{k}
\frac1{2\pi i}\int_{\calC _{\xi}}(z-z_0)^{-m-1}e^{tz}(-z)^l\notag\\
\label{eq6.14}
&{}\times\frac1{(k+l+1)!}\frac{(\frac{d}{dz})^{k+l+1}Q(z,\xi)}{Q(z,\xi)}\,dz\\\noalign{\goodbreak}
={}&\sum_{l=0}^{m-1-k}\binom{k+l}{k}
\frac1{2\pi i}\int_{\calC _{\xi}}(z-z_0)^{-m-1}e^{tz}(-z)^l\notag\\
\label{eq6.15}    
&{}\times\tau_{k+l+1}\bigg(\frac1{z-\lambda_1(\xi)},\ldots,\frac1{z-\lambda_m(\xi)}\bigg)
\,dz.  
\end{align}
Above, $\lambda_1(\xi),\ldots,\lambda_m(\xi)$ is the sequence of
eigenvalues of $\wid P(\xi)$  in which the number of occurrences of each eigenvalue is equal to its spectral multiplicity, and \vg
$$
Q(z,\xi)=\det(z\1-\wid P(\xi))=\prod_{k=1}^m(z-\lambda_k(\xi))
$$
is the characteristic polynomial of $\wid P(\xi)$.

By \eqref{eq6.13} the proposition follows once it is proved that
\begin{equation}\label{eq6.16}
a_k\in C^\infty(\bbR^{1+n};\bbC)
\end{equation}  
and
\begin{equation}\label{eq6.17}
\sup\{e^{-(\omega_0+\varepsilon)t}|a_k(t,\xi)|:t\in[0,\infty\[,\,
\xi\in\bbR^n\}<\infty
\end{equation}       
whenever $k=0,\ldots,m-1$ and $\varepsilon>0$.                       

In order to prove \eqref{eq6.16} notice that every $\xi_0\in\bbR^n$ has an open neighbourhood $U$  such that $\calC _{\xi_0}$  winds once about $\sigma(\wid P(\xi))$ whenever $\xi\in U$. (This follows from a theorem of Hurwitz. See \cite{S-Z}, Sec.~III.11.) Consequently, for 
every $\xi\in U$ one can replace $\calC _{\xi}$ by $\calC _{\xi_0}$
without changing the values of the integrals in \eqref{eq6.14}, and then \eqref{eq6.16} follows because $Q(z,\xi)$ is a $C^\infty$-function of $(z,\xi)$ non-vanishing on the open set
$\{(z,\xi)\in\bbC\times\bbR^n:z\not\in\sigma(\wid P(\xi)))\}$ which contains $\{(z,\xi)\in\bbC\times\bbR^n:z\in \calC _{\xi_0}\}$ if $\xi\in U$.

It remains to prove \eqref{eq6.17}. To this end, fix $\varepsilon>0$ and take $\delta\in\]0,\min(\varepsilon,\allowbreak\tfrac12(\hRe z_0
-\omega_0))]$. Let $\xi\in\bbR^n$. Since 
$\sigma(\wid P(\xi)))\subset\{z\in\bbC:\hRe z\le\omega_0\}$, without changing the values of the integrals in \eqref{eq6.15} one can choose a closed rectifiable path $\calC _{\xi}$ winding once about
$\sigma(\wid P(\xi))$ such that\vadjust{\vskip-2pt}
$$\vadjust{\vskip-2pt}
\calC _{\xi}\subset D_{\xi,\delta}:=\{z\in\bbC:\hRe z-\omega_0\le
\delta\le{\dist}(z,\sigma(\wid P(\xi)))\}.
$$
For every $\xi\in\bbR^n$ the straight line\vadjust{\vskip-2pt}
$$\vadjust{\vskip-2pt}
L=\{z\in\bbC:\hRe z=\omega_0+\delta\}
$$
is contained in $D_{\xi,\delta}$. Furthermore, whenever  $t\in[0,\infty\[$, $\xi\in\bbR^n$, $z\in D_{\xi,\delta}$ and $k+l=0,\ldots,m-1$, then
$|z-\lambda_i(\xi)|\ge\delta$, $|z-z_0|\ge\hRe z_0-\hRe z\ge\omega_0
+2\delta-(\omega_0+\delta)=\delta$ and\vadjust{\vskip-2pt}
$$
\bigg|\frac{z}{z-z_0}\bigg|=\bigg|1+\frac{z_0}{z-z_0}\bigg|
\le 1+\frac{|z_0|}{\delta},
$$
so that
$$\displaylines{\indent
|(z-z_0)^{-m-1}e^{tz}(-z)^l(z-\lambda_{i_1}(\xi))^{-1}\cdots
(z-\lambda_{i_{k+l+1}}(\xi))^{-1}|\hfill\cr
\hfill\eqalign{
&\le\bigg|\frac{z}{z-z_0}\bigg|^l|z-z_0|^{l+1-m}\delta^{-k-l-1}|z-z_0|^{-2}e^{(\omega_0+\delta)t}\cr
&\le \bigg(1+\frac{|z_0|}{\delta}\bigg)^l\delta^{-m-k}|z-z_0|^{-2}e^{(\omega_0+\delta)t},
\cr}\indent\cr}
$$
and hence
$$\displaylines{\indent
\bigg|(z-z_0)^{-m-1}e^{tz}(-z)^l\tau_{k+l+1}\bigg(\frac{1}{z-\lambda_1(\xi)},\ldots,\frac{1}{z-\lambda_m(\xi)}\bigg)\bigg|\hfill\cr
\hfill\le K_{\delta}|z-z_0|^{-2}e^{(\omega_0+\delta)t}\indent\cr}
$$
where $K_\delta$ is a finite number depending only on~$\delta$. Therefore, by Cauchy's integral theorem, the integration contour $\calC_\xi$ in \eqref{eq6.15} may be replaced by the straight line~$L$. Since $\delta\le\varepsilon$, one obtains the estimate
$$
|a_k(t,\xi)|\le \bigg(\sum_{l=0}^{m-1-k}\binom{k+l}{k}\bigg)
\frac{K_\delta}{2\pi}
\bigg(\int_L|z-z_0|^{-2}\,dz\bigg)e^{(\omega_0+\varepsilon)t}
$$
for $k=0,\ldots,m-1$, $t\in[0,\infty\[$ and $\xi\in\bbR^n$, proving \eqref{eq6.17}.      
\end{pf}

\section{Proof of Theorem 1 in cases (iv) and (v)}\label{sec7}

Theorem 1 is the conjunction of three implications: ${\rm(a)\Rightarrow(b)}\cap(\omega_E
\le \omega_0)$, 
(b)$\Rightarrow$(a), and (b)$\wedge(\omega_E
<\infty)\Rightarrow(\omega_0
\le \omega_E)$, the proofs of which differ in particular cases (i)--(v). The present section is devoted to Theorem 1(iv) and 1(v). The proofs in these cases are independent of the general theory of l.c.v.s.~and the advanced theory of distributions. 

Let either $E=(H^\infty(\bbR^n))^m$ and $\|u\|_j=(\sum_{0\le|\alpha|\le j}\|(\partial/\partial x)^\alpha u\|_{(L^2(\bbR^n)^m}^2)^{1/2}$ for  $j = 0, 1, \ldots,$
or $E=\{u\in(L^2(\bbR^n))^m:(P(\partial/\partial x))^lu\in(L^2(\bbR^n))^m$
for $l=1,2,\ldots\}$
and $\|u\|_j=\sum_{l=0}^j\|(P(\partial/\partial x))^lu\|_{(L^2(\bbR^n))^m}$
for  $j = 0, 1, \ldots .$ 

\medskip
\noindent
{\bf Proof of ${\rm(a)}\Rightarrow{\rm(b)}\wedge (\omega_E\le\omega_0)$.}
Suppose that (a) holds, i.e. $\omega_0\le\infty$. Let $(G_t)_{t\in\bbR}$ be the one-parameter group \eqref{eq2.1}. Whenever $E=(H^\infty(\bbR^n))^m$, $\varphi\in(C_c^\infty(\bbR^n))^m$, 
$u=\calF\varphi\in(Z_n)^m$, $j = 0, 1,\ldots$  and $t\in\bbR$, then, by the Remark at the end of Section~\ref{sec2}, 
$(\partial/\partial x)^{\alpha}G_tu=G_t(\partial/\partial x)^{\alpha} u\in(Z_n)^m$
and, by Plancherel's theorem, 
\begin{align*}
\|G_tu\|_j&=\bigg(\sum_{0\le|\alpha|\le j}\bigg\|G_t\bigg(\frac\partial{\partial x}\bigg)^\alpha u\bigg\|_0^2\bigg)^{1/2}\\
&=(2\pi)^{n/2}\bigg(\sum_{0\le|\alpha|\le j}\int_{\bbR^n}\|e^{t\wid P(\xi)}\xi^{\alpha}
\varphi(\xi)\|_{\bbC^m}^2\,d\xi\bigg)^{1/2}\\
&\le(2\pi)^{n/2}\bigg(\int_{\bbR^n}\|e^{t\wid P(\xi)}\|_{L\bbC^m}^2
\bigg(\sum_{0\le|\alpha|\le j}|\xi^{\alpha}|^2\bigg)\|\varphi(\xi)\|_{\bbC^m}^2\,d\xi\bigg)^{1/2}.
\end{align*}
Hence, by the Gelfand--Shilov inequality \eqref{eq6.9}, for every
$\varepsilon>0$ there are $K_\varepsilon\in\]0,\infty\[$ and $K'_\varepsilon\in\]0,\infty\[$ such that if $t\in[0,\infty\[$, $j=0,1,\ldots$ and $u\in(Z_n)^m$, then
$$\displaylines{
\|G_tu\|_j\hfill\cr
\hfill\eqalign{
\le{}&(2\pi)^{n/2} K_\varepsilon e^{(\omega_0+ \varepsilon)t}
\bigg(\!\int_{\bbR^n}(1+\|\wid P(\xi)\|^{m-1}_{L(\bbC^m)})^2 
\bigg(\sum_{0\le|\alpha|\le j}|\xi^{\alpha}|^2  \bigg) \|\varphi(\xi)\|_{\bbC^m}^2\,d\xi\!\bigg)^{1/2} \cr
\le{}&(2\pi)^{n/2} K'_\varepsilon  e^{(\omega_0+ \varepsilon)t} \bigg(\int_{\bbR^n}
\bigg(\sum_{0\le|\alpha|\le j+(m-1)d} |\xi^{\alpha}|^2\bigg)\|\varphi(\xi)\|_{\bbC^m}^2\,d\xi\bigg)^{1/2}\cr
\le{}&(2\pi)^{n/2} K'_\varepsilon  e^{(\omega_0+ \varepsilon)t} 
\bigg(\sum_{0\le|\alpha|\le j+(m-1)d}\bigg\|\calF^{-1}
\bigg(\frac\partial{\partial x}\bigg)^\alpha u\bigg\|_{(L^2(\bbR^n))^m}^2\bigg)^{1/2}.
\cr}\cr}
 $$
Again by Plancherel's theorem, it follows that
\begin{align}
\|G_tu\|_j&\le K'_\varepsilon  e^{(\omega_0+ \varepsilon)t} 
\bigg(\sum_{0\le|\alpha|\le j+(m-1)d}\bigg\|
\bigg(\frac{\partial}{\partial x}\bigg)^\alpha u\bigg\|_{(L^2(\bbR^n))^m}^2\bigg)^{1/2}\notag\\
\label{eq7.1}
&=K'_\varepsilon  e^{(\omega_0+ \varepsilon)t} \|u\|_{j+(m-1)d}
\end{align}
whenever $t\in[0,\infty\[$,  $u\in(Z_n)^m$ and $j=0,1,\ldots.$ Since
$(Z_n)^m$ is dense in  $E$, one concludes that $G_tE\subset E$ for every $t\in[0,\infty\[$, and that

\tekstl{the operators $S_t=G_t|_E$, $t\in[0,\infty\[$, 
constitute a one-parameter 
semigroup $(S_t)_{t\ge 0}\subset L(E)$ such that $\omega_E\le\omega_0$.}{eq7.2}

Now we are going to prove that \eqref{eq7.2} holds also in case (v), i.e. when
$E=\{u\in(L^2(\bbR^n))^m:(P(\partial/{\partial x}))^l u\in(L^2(\bbR^n))^m$
for $l=1,2,\ldots\}$. In this case for every $\varphi\in(C_c(\bbR^n))^m$, $u\in\calF
\varphi\in (Z_n)^m$, $j,l=0,1,\ldots$ and $t\in\bbR$ one has
$(P(\partial/{\partial x}))^l G_t u=G_t(P(\partial/{\partial x}))^l u\in (Z_n)^m$
and, by Plancherel's theorem,  
\begin{align*}
\|G_tu\|_j&=\sum_{l=0}^j\bigg\|\bigg(P\bigg(\frac\partial{\partial x}\bigg)\bigg)^l
G_t u\bigg\|_{(L^2(\bbR^n))^m}\\
&=\sum_{l=0}^j\bigg\|\bigg(P\bigg(\frac\partial{\partial x}\bigg)\bigg)^l\calF
e^{t\wid P}\varphi\bigg\|_{(L^2(\bbR^n))^m}\\
&=\sum_{l=0}^j\|\calF\wid P^le^{t\wid P}\varphi\|_{(L^2(\bbR^n))^m}
=(2\pi)^{n/2}\sum_{l=0}^j\|\wid P^l e^{t\wid P}\varphi\|_{(L^2(\bbR^n))^m}.
\end{align*}
Hence, by the Proposition at the end of Section~\ref{sec6}, for every  $\varepsilon>0$ there is $K_{\varepsilon}\in\]0,\infty\[$ such that whenever $\varepsilon\in
(C_c(\bbR^n))^m$, $u=\calF\varphi\in(Z_n)^m$, $j=0,1,\ldots$ and $t\in[0,\infty\[$, 
then\vadjust{\vskip-2pt}
\begin{align*}
\|G_tu\|_j&=(2\pi)^{n/2}\sum_{l=0}^j\bigg\|\sum_{k=0}^{2m}p_k(t,\bullet)\wid P^{k+l}
\varphi\bigg\|_{(L^2(\bbR^n))^m}\\[-2pt]
&\le(2\pi)^{n/2}K_{\varepsilon}e^{(\omega_0+ \varepsilon)t}\sum_{l=0}^{j+2m}\|\wid P^l\varphi\|_{(L^2(\bbR^n))^m}\\[-2pt]
&=K_{\varepsilon}e^{(\omega_0+ \varepsilon)t}\sum_{l=0}^{j+2m}(2\pi)^{n/2}
\bigg\|\calF^{-1}\bigg(P\bigg(\frac\partial{\partial x}\bigg)\bigg)^l u\bigg\|_{(L^2(\bbR^n))^m}.\vadjust{\vskip-2pt}
\end{align*}

Again by Plancherel's theorem, it follows that\vadjust{\vskip-2pt}
\begin{align}\label{eq7.3}
\|G_tu\|&\le K_{\varepsilon}e^{(\omega_0+ \varepsilon)t}\sum_{l=0}^{j+2m}
\bigg\|\bigg(P\bigg(\frac\partial{\partial x}\bigg)\bigg)^l u\bigg\|_{(L^2(\bbR^n))^m}\notag\\[-2pt]
&=K_\varepsilon e^{(\omega_0+\varepsilon)t}\|u\|_{j+2m}\vadjust{\vskip-2pt}
\end{align}  
whenever $t\in[0,\infty\[$, $u\in(Z_n)^m$ and $j=1,2,\ldots.$ Since, also in case (v),
$(Z_n)^m$ is dense in  $E$, one concludes that \eqref{eq7.2} holds.

It is easy to see that for both the Fr\'echet spaces $E$ considered in the present section,
$(C_c^\infty(\bbR^n))^m$ is continuously imbedded in $\calF^{-1}E$ and is sequentially dense in
$\calF^{-1}E$. Hence $(Z_n)^m=(\calF C_c^\infty(\bbR^n))^m$ is continuously imbedded
in $E$ and $(Z_n)^m$ is sequentially dense in~$E$. Therefore
in order to complete the proof of the implication ${\rm(a)\Rightarrow(b)}\wedge(\omega_E\le
\omega_0)$ it remains to apply to the one-parameter semigroup satisfying \eqref{eq7.2} the following

\medskip

\noindent{\bf Lemma.} {\it Suppose that {\rm\eqref{eq7.2}} holds for some l.c.v.s.~$E$  imbedded in $(Z'_n)^m$ such that $\frac\partial{\partial x_\nu}|_E\in L(E)$ for every $\nu=1,\ldots,n$. Suppose moreover that
$(Z_n)^m$ is continuously imbedded in $E$ and
sequentially dense in  $E$. Then $(S_t)_{t\ge0}\subset L(E)$
is a $(C_0)$-semigroup with infinitesimal generator $P(\partial/\partial x)|_E$.
}

\begin{pf}Pick $u\in E$ and let $(u_\nu)_{\nu=1,2,\ldots}\subset(Z_n)^m$ be a sequence converging to $u$ in the topology of~$E$. If $\nu\in\bbN$,
$t\in\]0,\infty\[$ and $\tau\in[0,\infty\[$, then  
$G_t u_\nu-G_{\tau} u_\nu=\int_{\tau}^t G_s P(\partial/\partial x)u_\nu\,ds$  and 
\begin{multline*}
\frac1t(G_t u_\nu-u_\nu)-P\bigg(\frac\partial{\partial x}\bigg)u_\nu
=\frac1t\int_0^t\bigg[G_sP\bigg(\frac\partial{\partial x}\bigg)u_\nu-P\bigg(\frac\partial{\partial x}\bigg)u_\nu\bigg]\,ds\\[-2pt]
=\frac1t\int_0^t\int_0^sG_\tau\bigg(P\bigg(\frac\partial{\partial x}\bigg)\bigg)^2u_\nu\,d\tau\,
ds
=\frac1t\int_0^t sG_{t-s}\bigg(P\bigg(\frac\partial{\partial x}\bigg)\bigg)^2u_\nu\,ds,
\vg\end{multline*}
 so that
$$
S_t u-S_\tau u=\lim_{\nu\to0}\int_\tau^t S_s P\bigg(\frac\partial{\partial x}\bigg)u_\nu\,ds
$$
and
$$
\frac1t(S_tu-u)-P\bigg(\frac\partial{\partial x}\bigg)u
=\lim_{\nu\to\infty}\frac1t\int_0^t s S_{t-s}\bigg(P\bigg(\frac\partial{\partial x}\bigg)\bigg)^2u_\nu\,ds.
$$
Here the integrals are Riemann integrals of continuous functions taking values in the complete l.c.v.s.~$(Z_n)^m$ (isomorphic to $(C_c^\infty(\bbR^n))^m$ which is complete; see [S2], p.~66, Theorem~I).  Let  $p$ be any continuous seminorm on~$E$. 
The restriction $p|_{(Z_n)^m}$ is a continuous seminorm on $(Z_n)^m$, so that
$$\displaylines{
p(S_tu-S_\tau u)=\lim_{\nu \to\infty}p
\bigg(\int_\tau^t S_sP\bigg(\frac\partial{\partial x}\bigg)u_\nu\,ds \bigg)
\le a|t-\tau|,\cr
a=\sup\{p(S_sP(\partial/\partial x)u_{\nu }):s\in[0,T],\,\nu \in\bbN\},\cr}
$$
and
$$\displaylines{
p\bigg(\frac1t(S_tu- u)-P\bigg(\frac\partial{\partial x}\bigg)u\bigg)
=\lim_{\nu \to\infty}p\bigg(\frac1t
\int_0^t s S_{t-s}\bigg(
P\bigg(\frac\partial{\partial x}\bigg)\bigg)u_\nu\,ds \bigg)
\le\frac12 bt,\cr
b=\sup\{p(S_s(P(\partial/\partial x))^2u_{\nu }):s\in[0,T],\,\nu \in\bbN\}.\cr}
$$ 
By \eqref{eq7.2} the semigroup $(S_t)_{t\ge0}\subset L(E)$ is locally equicontinuous, whence $\{S_s(P(\partial/\partial x))^ku_{\nu }):s\in[0,T],\,\nu \in\bbN\}$,  $k=1,2,\ldots,$ are bounded subsets of~$E$, and so  $a$  and  $b$  are finite. Consequently, the function $[0,\infty\[\ni t\mapsto S_tu\allowbreak\in E$ is continuous, and $\lim_{t\downarrow 0}
\frac1t(S_tu-u)=P(\partial/\partial x)u$
in the topology of~$E$. This proves that $(S_t)_{t\ge0}\subset L(E)$ is a $(C_0)$-semigroup with infinitesimal generator $P(\partial/\partial x)|_E$. 
\end{pf}

\noindent{\bf Proof of  ${\rm(b){\Rightarrow}(a)}$ and {\rm(b)}$\wedge(\omega_E<\infty)
\Rightarrow(\omega_0\le \omega_E)$.} If (b) holds, then, by the Uniqueness Theorem from Section~\ref{sec2}, for every $t\in[0,\infty\[$ there are $C(t)\in\]0,\infty\[$ and
$j(t)\in\bbN_0$ such that $\|G_tu\|_0=\|S_tu\|_0\le C(t)\|u\|_{j(t)}$ for every $u\in E$. If $\varphi\in(C_c^\infty(\bbR^n))^m$ and $u=\calF\varphi\in(Z_n)^m$, then 
$G_tu=\calF e^{t\wid P}\varphi$, so that
$$\eqalign{
\|e^{t\wid P}\varphi\|_0
&=\|\calF^{-1}G_tu\|_0\le(2\pi)^{-n/2}\|G_t u\|_0\le(2\pi)^{-n/2}
C(t)\|u\|_{j(t)}\cr
&=(2\pi)^{-n/2}C(t)\|\calF\varphi\|_{j(t)}
\le\bigg(\int_{\bbR^n}Q_t(\eta)^2\|\varphi(\eta)\|_{\bbC^m}^2\,d\eta\bigg)^{1/2}\cr}
$$ 
 where $Q_t$ is a real polynomial on~$\bbR^n$. Consequently, for every
$t\in[0,\infty\[$ there are $K(t)\in\]0,\infty\[$ and $k(t)\in\bbN$ such that
\begin{equation}\label{eq7.4}
\bigg(\int_{\bbR^n}\|e^{t\wid P}\varphi(\eta)\|_{\bbC^m}^2\,d\eta\bigg)^{1/2}
\le K(t)\bigg(\int_{\bbR^n}(1+|\eta|^2)^{k(t)}\|\varphi(\eta)\|_{\bbC^m}^2\,d\eta\bigg)^{1/2}
\end{equation}  
whenever   $\varphi\in(C_c^\infty(\bbR^n))^m$. Furthermore, if
${\rm(b)}\wedge(\omega_E<\infty)$ holds, then for every $\varepsilon>0$ there are
$K_\varepsilon\in\]0,\infty\[$ and $j_\varepsilon\in\bbN$ such that $(2\pi)^{-n/2}C(t)\le K_\varepsilon e^{(\omega_E+\varepsilon)t}$ and $j(t)\le j_\varepsilon$
for every $t\in[0,\infty\[$, whence it follows that for some  $k_\varepsilon\in\bbN$ one has
$$\displaylines{
\bigg(\int_{\bbR^n}\|e^{t\wid P}\varphi(\eta)\|_{\bbC^m}^2\,d\eta\bigg)^{1/2}
\le K_\varepsilon e^{(\omega_E+\varepsilon)t}
\bigg(\int_{\bbR^n}(1+|\eta|^2)^{k_\varepsilon}
\|\varphi(\eta)\|_{\bbC^m}^2\,d\eta\bigg)^{1/2}\hfill
\refstepcounter{equation}(\theequation)\label{eq7.5}\cr}
$$
whenever $t\in[0,\infty\[$ and $\varphi\in(C_c^\infty(\bbR^n))^m$. Inequalites \eqref{eq7.4} and \eqref{eq7.5} resemble \eqref{eq4.4}, and our subsequent arguments resemble that from the proof of Theorem~2 of Section~4. All this is similar to the argument used by
T.~Ushijima in the Correction to~\cite{U}.

For every $t\in[0,\infty\[$ and $\xi\in\bbR^n $ choose $z_{t,\xi}\in\bbC^m$ such that
$\|z_{t,\xi}\|_{\bbC^m}=1$ and $\|e^{t\wid P(\xi)}\|_{L(\bbC^m)}=
\|e^{t\wid P(\xi)}z_{t,\xi}\|_{\bbC^m}$. Let
$(\phi_{\xi,\nu })_{\nu =1,2,\ldots}\subset C_c^\infty(\bbR^n)$ be a sequence of non-negative functions such that for every
$\nu =1,2,\ldots$ the support of $\phi_{\xi,\nu }$ is contained in the ball with centre
 at $\xi$ and radius $1/\nu $, and $\int_{\bbR^n}\phi_{\xi,\nu }(\eta)^2\,d\eta=1$.

If (b) holds, then applying \eqref{eq7.4} to $\varphi(\eta)
=\phi_{\xi,\nu }(\eta)z_{t,\xi}$ one concludes that whenever  $t\in[0,\infty\[$ and $\xi\in\bbR^n$, then
\begin{align}\label{eq7.6}
\|e^{t\wid P(\xi)}\|_{L(\bbC^m)}=\|e^{t\wid P(\xi)}z_{t,\xi}\|_{\bbC^m}
&=\lim_{\nu \to\infty}\bigg(\int_{\bbR^n}\|e^{t\wid P(\eta)}\phi_{\xi,\nu }(\eta)z_{t,\xi}\|_{\bbC^m}^2
\,d\eta\bigg)^{1/2}\notag\\
&\le K(t)\lim_{\nu \to\infty}\bigg(\int_{\bbR^n}(1+|\eta|^2)^{k(t)}\phi_{\xi,\nu }(\eta)^2
\,d\eta\bigg)^{1/2}\notag\\
&= K(t)(1+|\xi|^2)^{k(t)/2}.
\end{align}      
Let $\rho(e^{\wid P(\xi)})$ denote the spectral radius of the matrix $e^{\wid P(\xi)}$. By Corollary 2.4 on 
p.~252 of \cite{E-N}, \eqref{eq7.6} implies that
$$\eqalign{
\max\{\hRe \lambda:\lambda\in\sigma(\wid P(\xi))\}
&=\log\rho(e^{\wid P(\xi)})
\le\log\|e^{\wid P(\xi)}\|_{L(\bbC^m)}\cr
&\le\log K(1)+\frac12k(1)\log (1+|\xi|^2)\cr}
$$                                                            
for every $\xi\in\bbR^n$. By the Corollary at the end of our Section~\ref{sec5}, it follows that
$$
\omega_0=\sup\{\hRe\lambda:\lambda\in\sigma(\wid P(\xi)),\,\xi\in\bbR^n\}<\infty.
$$

If ${\rm(b)}\wedge(\omega_E<\omega)$ holds, then the difficult results quoted in 
Section~\ref{sec5} need not be used. Applying \eqref{eq7.5} to
$\varphi(\eta)=\phi_{\xi,\nu }(\eta)z_{t,\xi}$ one concludes that whenever
$t\in[0,\infty\[$ and $\xi\in\bbR^n$, then
$$
\|e^{t\wid P(\xi)}\|_{L(\bbC^m)}\le K_\varepsilon e^{(\omega_\varepsilon+\varepsilon)t}
(1+|\xi|^2)^{k_\varepsilon/2},
$$
whence, by Proposition 2.2, p.~251, and Corollary 2.4, p.~252, in \cite{E-N},
$$
\max\{\hRe \lambda:\lambda\in\sigma(\wid P(\xi))\}
=\lim_{t\to\infty}\frac1t\log\|e^{t\wid P(\xi)}\|_{L(\bbC^m)}
\le\omega_E+\varepsilon.
$$
This implies that $\omega_0\le\omega_E$ because $\varepsilon>0$ is arbitrary.

\section{Conditions on  $e^{t\wid P(\xi)}$ equivalent to the Petrovski\u\i\  correctness}\label{sec8}

Let $\wid P(\xi)$ be the symbol of the matricial differential operator $P(\partial/\partial x)$
defined in Section~\ref{sec1}. For any $\omega\in\bbR$ consider the conditions:
$$\displaylines{\indent
\sup\{\hRe\lambda:\lambda\in\sigma(\wid P(\xi)),\,\xi\in\bbR^n\}\le\omega\mbox{ (the Petrovski\u\i\  correctness)};\hfill
\refstepcounter{equation}(\theequation)\label{eq8.1}\cr
\noalign{\vskip6pt}
\indent
\mbox{there is $k\in\bbN$ such that}\hfill\cr
\indent
\sup\{e^{-(\omega+\varepsilon)t}(1+|\xi|)^{-k}\|e^{t\wid P(\xi)}\|_{M_m}:
t\in[0,\infty\[,\,\xi\in\bbR^n\}<\infty\hfill
\refstepcounter{equation}(\theequation)\label{eq8.2}\cr
\indent
\mbox{for every $\varepsilon>0$;}\hfill\cr\noalign{\goodbreak}
\noalign{\vskip6pt}
\indent
\mbox{for every multiindex $\alpha\in\bbN_0^n$ there is $k_{\alpha}\in\bbN$ such that}\hfill\cr
\indent
\sup\{e^{-(\omega+\varepsilon)t}(1+|\xi|)^{-k_\alpha}\|(\partial/\partial\xi)^\alpha
e^{t\wid P(\xi)}\|_{M_m}:t\in[0,\infty\[,\,\xi\in\bbR^n\}<\infty\hfill
\refstepcounter{equation}(\theequation)\label{eq8.3}\cr
\indent
\mbox{for every  $\varepsilon>0$.}\hfill\cr}
$$
Condition \eqref{eq8.1} implies \eqref{eq8.2} by the Gelfand--Shilov inequality \eqref{eq6.9}, and \eqref{eq8.2} implies \eqref{eq8.1} because
$\max\{\hRe \lambda:\lambda\in\sigma(\wid P(\xi))\}=t^{-1}\log \rho(e^{t\wid P(\xi)})
\le t^{-1}\log \|e^{t\wid P(\xi)}\|_{L(\bbC^m)}$ for every
$t\in\]0,\infty\[$ where $\rho$ stands for the spectral radius. See \cite{E-N}, p.~252. 
The partial derivatives occurring in \eqref{eq8.3} make sense because the function
$\bbR^{1+n}\ni(t,\xi)\mapsto e^{t\wid P(\xi)}\in M_m$
is infinitely differentiable, by arguments mentioned at the beginning of Section~\ref{sec2}. Obviously \eqref{eq8.3} implies \eqref{eq8.2}, and the proof of the converse implication will be given shortly. Therefore for any fixed $\omega\in\bbR$ the conditions \eqref{eq8.1}, \eqref{eq8.2} and \eqref{eq8.3} are equivalent.

I. G. Petrovski\u\i\  considered in [P] the following conditions:
$$\displaylines{\indent
\sup\{(1+\log(1+|\xi|))^{-1}\hRe\lambda:\lambda\in\sigma(\wid P(\xi)),\,\xi\in\bbR^n\}<\infty;\hfill
\refstepcounter{equation}(\theequation)\label{eq8.4}\cr
\noalign{\vskip6pt}
\indent
\mbox{for every $T\in\]0,\infty\[$ there is $k_T\in\bbN$ such that}\notag\hfill\cr
\indent
\sup\{(1+|\xi|)^{-k_T}\|e^{t\wid P(\xi)}\|_{M_m}:t\in[0,T],\,\xi\in\bbR^n\}<\infty;\hfill
\refstepcounter{equation}(\theequation)\label{eq8.5}\cr
\noalign{\vskip6pt}
\indent
\mbox{for every multiindex $\alpha\in\bbN_0^n$ and every $T\in\]0,\infty\[$}\hfill\cr
\indent
\mbox{there is $k_{\alpha,T}\in\bbN$ such that}\hfill\cr
\indent
\sup\{(1+|\xi|)^{-k_{\alpha,T}}\|(\partial/\partial\xi)^\alpha
e^{t\wid P(\xi)}\|_{M_m}:t\in[0,T],\,\xi\in\bbR^n\}<\infty.\hfill
\refstepcounter{equation}(\theequation)\label{eq8.6}\cr}
$$  
The three conditions \eqref{eq8.4}--\eqref{eq8.6} are equivalent to 
each other, and each is equivalent to the existence of an $\omega\in\bbR$ for which the conditions \eqref{eq8.1}--\eqref{eq8.3} are satisfied. This follows from the Corollary at the end of Section 5 and arguments similar to those proving the mutual equivalence of \eqref{eq8.1}, \eqref{eq8.2} and \eqref{eq8.3}.

\medskip

\noindent
{\bf Proof of the implication {\rm\eqref{eq8.2}}$\Rightarrow${\rm\eqref{eq8.3}}.} For every $\alpha=(\alpha_1,\ldots,\alpha_n)\in\bbN_0^n$, $\xi\in\bbR^n$ and $t\in[0,\infty\[$ put
$$\displaylines{
\wid P_{\alpha}(\xi)=\bigg(\frac\partial{\partial\xi_1}\bigg)^{\alpha_1}\cdots\bigg(\frac\partial{\partial\xi_n}\bigg)^{\alpha_n}\wid P(\xi),\cr
U_{\alpha}(t,\xi)=\bigg(\frac\partial{\partial\xi_1}\bigg)^{\alpha_1}\cdots\bigg(\frac\partial{\partial\xi_n}\bigg)^{\alpha_n}e^{t\wid P(\xi)}.
\cr}
$$
If $\alpha,\beta\in\bbN_0^n$, then let $\beta\le\alpha$ mean that 
$\beta_\nu \le\alpha_\nu $ for every $\nu =1,\ldots,n$. If $\beta\le\alpha$, then 
$\binom{\alpha}{\beta}:=\binom{\alpha_1}{\beta_1}\cdots\binom{\alpha_n}{\beta_n}$ where 
$\binom{\alpha_\nu }{\beta_\nu }=\frac{\alpha_\nu !}{\beta_\nu !(\alpha_\nu -\beta_\nu )!}$. Condition \eqref{eq8.3} means that whenever $\alpha\in\bbN_0^n$ then
$$\displaylines{\indent
\mbox{there is $k\in\bbN$ such that}\hfill\cr
\indent
\sup\{e^{-(\omega+\varepsilon)t}(1+|\xi|)^{-k}\|U_\alpha(t,\xi)\|:
t\in[0,\infty\[,\,\xi\in\bbR^n\}<\infty\hfill(8.7)_\alpha
\refstepcounter{equation}
\label{eq8.7}\cr
\indent
\mbox{for every $\varepsilon>0$.}\hfill\cr}
$$
Condition \eqref{eq8.2} is identical with \eqref{eq8.7}$_0$. Hence the implication \eqref{eq8.2}$\Rightarrow$\eqref{eq8.3} will follow once we prove that if $l\in\bbN_0$ and \eqref{eq8.7}$_\beta$  holds for every $\beta\in\bbN_0^n$ such that $|\beta|=
\beta_1+\cdots+\beta_n\le l$, then \eqref{eq8.7}$_\alpha$  holds for every
$\alpha\in\bbN_0^n$ such that $|\alpha|=l+1$. So, pick any $\alpha$ such that $|\alpha|=l+1$. Then
\begin{align}\label{eq8.8}
\frac d{dt}U_\alpha(t,\xi)&=\bigg(\frac\partial{\partial\xi}\bigg)^{\alpha}
\frac d{dt}U_0(t,\xi)=\bigg(\frac\partial{\partial\xi}\bigg)^{\alpha}
(\wid P(\xi)U_0(t,\xi))\notag\\
&
=\wid P(\xi)U_\alpha(t,\xi)+
V_\alpha(t,\xi)
\end{align} 
where
$$
V_\alpha(t,\xi)=\sum_{\beta\le\alpha,|\beta|\le l}\binom{\alpha}{\beta}
\wid P_{\alpha-\beta}(\xi)
U_\beta(t,\xi).
$$
Since \eqref{eq8.7}$_\beta$ holds whenever $|\beta|\le l$, it follows that 
$$\displaylines{\indent
\mbox{there is $k\in\bbN$ such that}\hfill\cr
\indent
\sup\{e^{-(\omega+\varepsilon)t}(1+|\xi|)^{-k}\|V_\alpha(t,\xi)\|:
t\in[0,\infty\[,\,\xi\in\bbR^n\}<\infty\hfill	
\refstepcounter{equation}(\theequation)\label{eq8.9}\cr
\indent
\mbox{for every $\varepsilon>0$.}\hfill\cr}
$$
By \eqref{eq8.8} one has 
\begin{equation}\label{eq8.10}
U_\alpha(t,\xi)=\int_0^t U_0(t-\tau,\xi)V_\alpha(\tau,\xi)\,d\tau,\ \quad
t\in[0,\infty\[,\,\xi\in\bbR^n.
\end{equation}
Conditions \eqref{eq8.7}$_0$  and \eqref{eq8.9} imply \eqref{eq8.7}$_\alpha$, by \eqref{eq8.10}. 

\medskip

\noindent
{\bf Remark.} The above proof is similar to the proof of Lemma~2 in Sec.~2 of Chap.~1 of \cite{P}. By \eqref{eq6.9}, condition \eqref{eq8.1} implies \eqref{eq8.2} with $k=(m-1)d$. From this last, by the induction procedure used above, one obtains \eqref{eq8.3} with 
$k_\alpha=(md-1)(|\alpha|+1)$.

\section{The space $\calO_{M}(\bbR^n;M_m)$}\label{sec9}

The paper \cite{P} of I.~G.~Petrovski\u\i\  makes evident the fundamental role of smooth slowly increasing functions in the theory of Cauchy's problem for systems of PDE with constant coefficients. A continuous function $\phi$ defined on $\bbR^n$ is called
{\it slowly increasing}\/ if there is $k\in\bbN_0$ such that $\sup\{(1+|\xi|)^{-k}
|\phi(\xi)|:\xi\in\bbR^n\}<\infty$. The space 
$\calO_M(\bbR^n;M_m)$ of $M_m$-valued slowly increasing infinitely differentiable functions on $\bbR^n$  consists of $M_m$-valued $C^\infty$-functions $\phi$ on $\bbR^n$ such that $\phi$ and all its partial derivatives are slowly increasing. The present section is devoted to the properties of $\calO_M(\bbR^n;M_m)$ as the space of multipliers for
$(\calS(\bbR^n))^m$ and $(\calS'(\bbR^n))^m$. 

An $M_m$-valued function $\phi$ defined on $\bbR^n$ is called a {\it multiplier for} $(\calS(\bbR^n))^m$ if $\phi\bullet\varphi\in(\calS(\bbR^n))^m$ for every
$\varphi\in(\calS(\bbR^n))^m$ and the multiplication operator
$\phi\,\bullet:\varphi\mapsto\phi\cdot\varphi$ belongs to $L((\calS(\bbR^n))^m)$. The set of multipliers for $(\calS(\bbR^n))^m$ will be denoted by $\calM((\calS(\bbR^n))^m)$. Obviously  $\calM((\calS(\bbR^n))^m)\subset C^\infty(\bbR^n;M_m)$.

Let $\calS'(\bbR^n)$ be the L.~Schwartz space of tempered distributions on~$\bbR^n$, and let $\calS'_b(\bbR^n)$ denote $\calS'(\bbR^n)$ endowed with topology of uniform convergence on bounded subsets of  $\calS(\bbR^n)$. A function
$\phi\in C^\infty(\bbR^{n};M_m)$ is called a
{\it multiplier for} $(\calS'_b(\bbR^n))^m$ if $\phi\cdot T
\in(\calS'(\bbR^n))^m$ for every $T=(T_1,\ldots,T_m)\in(\calS'(\bbR^n))^m$ and the multiplication operator
$\phi\,\bullet:T\mapsto \phi\cdot T$ belongs to $L((\calS'_b(\bbR^n))^m)$. The set of multipliers for $(\calS'_b(\bbR^n))^m$ will be denoted by $\calM((\calS'_b(\bbR^n))^m)$.

\medskip
\noindent
{\bf Proposition.} $\calM((\calS'_b(\bbR^n))^m)=\calM((\calS(\bbR^n))^m)$.

\begin{pf}If $\phi(t)=(\phi_{i,j}(t))_{i,j=1,\ldots,m}$, then
$\phi\in\calM((\calS(\bbR^n))^m)$ if and only if $\phi_{i,j}\in
\calM(\calS(\bbR^n))$ for every $i,j=1,\ldots,m$. A similar equivalence holds for $\calS'_b(\bbR^n)$.
Therefore it is sufficient to prove the Proposition for $m=1$.
It is clear that $\calM(\calS(\bbR^n))\subset\calM(\calS'_b(\bbR^n))$. 
To prove the opposite inclusion notice that the pair of 
l.c.v.s.~$(\calS(\bbR^n),\calS'_b(\bbR^n))$
is reflexive with respect to the duality form 
$\calS(\bbR^n)\times\calS'_b(\bbR^n)\ni(u,T)\mapsto
T(u)\in\bbC$. See \cite{S}, Sec.~VII.4, p.~238, remarks after Theorem IV;  \cite{Y}, p.~140, Theorem~2.

Henceforth fix any $\phi\in\calM(\calS'(\bbR^n))$. We have to prove that
 $\phi\in\calM(\calS(\bbR^n))$. Whenever $u\in\calS(\bbR^n)$, then
$\calS'_b(\bbR^n)\ni T\mapsto\<\phi\cdot T,u>\allowbreak\in\bbC$ is a linear functional continuous on $\calS'_b(\bbR^n)$ and therefore, by reflexivity, there is a unique
$v\in\calS(\bbR^n)$ such that $\<\phi\cdot T,u>=\<T,v>$ for every 
$T\in\calS'(\bbR^n)$. The map $A:\calS(\bbR^n)\ni u\mapsto
v\in\calS(\bbR^n)$ is algebraic linear. Let
$(u_k)_{k=1,2,\ldots}\subset\calS(\bbR^n)$ be a sequence  converging to zero in the Fr\'echet topology in $\calS(\bbR^n)$, and let $B$ be a bounded subset of $\calS'_b(\bbR^n)$. Then 
$\phi\bullet B$ is again bounded, so that
$$
\lim_{k\to\infty}\sup\{|\<T,Au_k>|:T\in B\}
=\lim_{k\to\infty}\sup\{|\<\phi\cdot T,u_k>|:T\in B\}=0.
$$
This means that $\lim_{k\to\infty}Au_k=0$ in the topology of
the strong dual of $\calS'(\bbR^n)$, and hence, again by reflexivity,
$\lim_{k\to\infty}Au_k=0$ in the 
 original Fr\'echet topology of $\calS(\bbR^n)$. Thus
\begin{equation}\label{eq9.1}
\<\phi\cdot T,u>=\<T,Au>\ \quad\mbox{for every }
u\in\calS(\bbR^n)
\ \mbox{and}\  
T\in\calS'(\bbR^n)
\end{equation}
where $A=L(\calS(\bbR^n)).$ Fix any $u\in\calS(\bbR^n)$. The Proposition will follow once we show that $\phi\bullet u=Au$. To prove this, let
$(u_k)_{k=1,2,\ldots}\subset C^\infty_c(\bbR^n)$ be a sequence that converges to  $u$  in the Fr\'echet topology of $\calS(\bbR^n)$. 
By the definition of multiplication of a distribution by a $C^\infty$-function, and 
by~\eqref{eq9.1}, 
\begin{multline*}
\lim_{k\to\infty}\sup\{|\<T,\phi\bullet u_k-Au>|:T\in B\}\\
\begin{aligned}
&=\lim_{k\to\infty}\sup\{|\<\phi\bullet T, u_k>-\<T,Au>|:T\in B\}\\
&=\lim_{k\to\infty}\sup\{|\<T,Au_k-Au>|:T\in B\}=0
\end {aligned} 
\end{multline*}                                                       for every bounded subset $B$ of $\calS'_b(\bbR^n)$. By reflexivity, this implies that
$\lim_{k\to\infty}\phi\bullet u_k=Au\in\calS(\bbR^n)$ in the Fr\'echet topology of
$\calS(\bbR^n).$ On the other hand, $\lim_{k\to\infty}\phi\bullet u_k=\phi\bullet u$ pointwise on~$\bbR^n$. Therefore $\phi\bullet u=Au\in\calS(\bbR^n)$ for every 
$u\in\calS(\bbR^n),$ 
which means that $\phi\in\calM(\calS(\bbR^n))$.

We now formulate the main results of the present section.
\end{pf}

\noindent
{\bf Theorem A.}  $\calO_M(\bbR^n;M_m)=\calM((\calS(\bbR^n))^m)$.

\medskip

This Theorem and the preceding Proposition imply at once 

\medskip
\noindent
{\bf Corollary.} $\calO_M(\bbR^n;M_m)=\calM((\calS'(\bbR^n))^m)$.
\medskip

This Corollary (for $m = 1$)  is mentioned without proof on p.~246 of \cite{S}. The Corollary and Theorem A (for $m = 1$) are formulated simultaneously in Theorem 25.5 stated without proof on p.~275 of~\cite{T}.

\medskip
\noindent
{\bf Theorem B.} {\it  For any subset  $B$ of $\calO_M(\bbR^n;M_m)$
the following five conditions are equivalent:
$$\displaylines{\indent
\mbox{for every $\alpha\in\bbN^n_0$ there is $k_\alpha\in\bbN^n_0$ such that}\hfill\cr
\indent
\sup\{(1+|\xi|)^{-k}\|(\partial/{\partial\xi})^{\alpha}\phi(\xi)\|_{M_m}:
\xi\in\bbR^n,\, \phi\in B\}<\infty;\hfill
\refstepcounter{equation}(\theequation)\label{eq9.2}\cr
\noalign{\vskip6pt}
\indent
\mbox{the family of multiplication operators}\hfill\cr
\indent
\{\phi\,\bullet :\phi\in B\}\subset L((\calS(\bbR^n))^m)
\ \quad\mbox{is equicontinuous;}\hfill
\refstepcounter{equation}(\theequation)\label{eq9.3}\cr
\noalign{\vskip6pt}
\indent
B\bullet C\mbox{ is a bounded subset of }(\calS(\bbR^n))^m\hfill\cr
\indent
\mbox{whenever $C$  is a bounded subset of $(\calS(\bbR^n))^m$};\hfill
\refstepcounter{equation}(\theequation)\label{eq9.4}\cr\noalign{\goodbreak}
\noalign{\vskip6pt}
\indent
\mbox{the family of multiplication operators}\hfill\cr 
\indent
\{\phi\,\bullet :\phi\in B\}\subset L((\calS'(\bbR^n))^m)
\mbox{ is equicontinuous;}\hfill
\refstepcounter{equation}(\theequation)\label{eq9.5}\cr
\noalign{\vskip6pt}
\indent
B\bullet C'\mbox{ is a bounded subset of }(\calS'(\bbR^n))^m\hfill\cr
\indent
\mbox{whenever $C'$ is a bounded subset of $(\calS'(\bbR^n))^m$.}\hfill
\refstepcounter{equation}(\theequation)\label{eq9.6}\cr}
$$\vskip-\lastskip
}\medskip

In Theorem A the inclusion  $\calO_M(\bbR^n;M_m)\subset\calM((\calS(\bbR^n))^m)$ is obvious and the difficult opposite inclusion is a particular case of the implication 
\eqref{eq9.3}$\Rightarrow$\eqref{eq9.2}. The implication \eqref{eq9.2}$\Rightarrow$\eqref{eq9.3} is obvious. For
$q=2,\ldots,6$ denote by $(9.q)^\dag$  the condition obtained from $(9.q)$ by replacing 
$B$  by $B^\dag=\{\phi^\dag:\phi\in B\}$ where $\phi^\dag$ denotes the transpose of the matrix~$\phi$. It is obvious that $\eqref{eq9.2}^\dag{\Leftrightarrow}\eqref{eq9.2}$. Furthermore, 
$(9.3){\Rightarrow}(9.4){\Rightarrow}(9.5)^\dag{\Rightarrow}(9.6)^\dag{\Rightarrow}(9.3)$. Indeed, the implications $\eqref{eq9.3}{\Rightarrow}\eqref{eq9.4}$ and
$\eqref{eq9.5}^\dag{\Rightarrow} \eqref{eq9.6}^\dag$ are trivial. The implication 
$(9.4){\Rightarrow}(9.5)^\dag$ follows from the fact that the polars of bounded subsets of
$(\calS(\bbR^n))^m$ constitute a basis of the topology of $(\calS'_b(\bbR^n))^m$. The implication $(9.6)^\dag{\Rightarrow}(9.3)$ follows from the fact  that
$\{(\calS(\bbR^n))^m,(\calS'_b(\bbR^n))^m\}$ is a reflexive pair of  l.c.v.s.~in duality, so that the polars of bounded subsets of $(\calS'(\bbR^n))^m$ constitute a basis for the original Fr\'echet topology of  $(\calS(\bbR^n))^m$. Hence  the equivalence 
$(9.2){\Rightarrow} (9.3)$ implies that $(9.q){\Rightarrow}(9.q)^\dag$  for $q=2,\ldots,6$. Therefore all what remains to do is to prove the implication 
$\eqref{eq9.3}{\Rightarrow} \eqref{eq9.2}$.

In order to simplify the proof of $\eqref{eq9.3}{\Rightarrow} \eqref{eq9.2}$, notice the
following equality concerning elements of $L((\calS(\bbR^n))^m)$:
$$
\bigg(\bigg(\frac\partial{\partial\xi}\bigg)^\alpha
\phi\bigg)\bullet=\bigg(\frac\partial{\partial\xi}\bigg)^{\alpha}(\phi\,\bullet)-\sum_{\beta\le\alpha,\,|\beta|<|\alpha|}\binom{\alpha}{\beta}
\bigg(\bigg(\bigg(\frac\partial{\partial\xi}\bigg)^{\beta}\phi\bigg)\bullet\bigg)
\bigg(\frac\partial{\partial\xi}\bigg)^{\alpha-\beta}.
$$
 By induction on
$|\alpha|=\alpha_1+\cdots+ \alpha_n$ this equality implies that if 
$\phi\in\calM((\calS(\bbR^n))^m)$, then
$(\partial/{\partial\xi})^{\alpha}\phi\in\calM((\calS(\bbR^n))^m)$ for every $\alpha\in\nobreak\bbN^n_0$, and if \eqref{eq9.3} holds then for every $\alpha\in\bbN^n_0$
the family of multiplication operators 
$\{((\partial/{\partial\xi})^{\alpha}\phi)\,\bullet:\phi\in B\}\subset L((\calS(\bbR^n))^m)$
is equicontinuous. Therefore the implication 
$\eqref{eq9.3}{\Rightarrow}\eqref{eq9.2}$ will follow once we prove that \eqref{eq9.3} implies the condition
$$\displaylines{\indent
\mbox{there is $k\in\bbN_0$ such that}\hfill\cr
\indent
\sup\{(1+|\xi|)^{-k}\|\phi(\xi)\|_{M_m}:\xi\in\bbR^n,\,\phi\in B\}<\infty.\hfill
\refstepcounter{equation}(\theequation)\label{eq9.7}
\cr}
$$
In the proof of (9.3)$\Rightarrow$(9.7) we will use

\medskip
\noindent 
{\bf Lemma.} $\calM((\calS(\bbR^n))^m)\subset\calS'(\bbR^n;M_m)$.

\begin{pf}Let $\phi\in\calM((\calS(\bbR^n))^m)$. Then
$\calF^{-1}(\phi\,\bullet)\calF\in L((\calS(\bbR^n))^m)$ and $\calF^{-1}(\phi\,\bullet)\calF$
commutes with translations. Therefore, by a variant of a theorem of L.~Schwartz (\cite{S}, p.~162; \cite{Y}, p.~158) there is a unique $T\in\calS'(\bbR^n;M_m)$ such that 
$\calF^{-1}(\phi\bullet\calF\varphi)=T*\varphi$ for every $\varphi\in(\calS(\bbR^n))^m$. 
Consequently, $\phi\bullet\calF\varphi=\calF(T*\varphi)=\calF T\bullet\calF\varphi$ for every $\varphi\in(\calS(\bbR^n))^m$, and hence 
$$
\phi=\calF T\in \calS'(\bbR^n;M_m).
$$
\end{pf}

\noindent{\bf Proof of ${\rm\eqref{eq9.3}{\Rightarrow}\eqref{eq9.7}}$.}
 This proof, or rather the proof of a Fourier precursor of 
${\rm\eqref{eq9.8}{\Rightarrow}\eqref{eq9.9}}$, resembles that of  Theorem 3.1 in \cite{Ch}, pp.~82--83, and that of a part of Theorem XXV in Sec.~VI.8 of~\cite{S}.

For notational convenience we introduce a set  $J$  of indices and a one-to-one map of  
$J\ni\iota\mapsto\phi_\iota\in B$ of $J$ onto~$B$.
We assume that $\phi_\iota\in\calM((\calS(\bbR^n))^m)$ for every $\iota\in J$ and that 
$$\displaylines{\indent	
\mbox{the family of multiplication operators}\hfill\cr
\indent	
\{\phi_\iota\,\bullet: \iota\in J\}\subset L((\calS(\bbR^n))^m)\mbox{ is equicontinuous.}
\hfill
\refstepcounter{equation}(\theequation)\label{eq9.8}
\cr}
$$ 
We have to prove that
$$\displaylines{\indent	
\mbox{there is $k\in\bbN_0$ such that}\hfill\cr
\indent
\sup\{(1+|\xi|)^{-k}\|\phi_\iota(\xi)\|_{M_m}:\xi\in\bbR^n,\,\iota\in J\}<\infty.\hfill
\refstepcounter{equation}(\theequation)\label{eq9.9}
\cr}
$$                             
This last condition will follow once we prove that there is a scalar polynomial  $Q$
and for every $\iota\in J$  there are $f_\iota,g_\iota\in L^1(\bbR^n;M_m)$ such that 
\begin{equation}\label{eq9.10}
\sup_{\iota\in J}\|f_\iota\|_{L^1(\bbR^n;M_m)}<\infty,\ \quad
\sup_{\iota\in J}\|g_\iota\|_{L^1(\bbR^n;M_m)}<\infty
\end{equation}
and
\begin{equation}\label{eq9.11}
T_\iota:=\calF\phi_\iota=Q  \bigg(\frac\partial{\partial x}\bigg) f_\iota+ g_\iota
\ \quad\mbox{for every }\iota\in J
\end{equation}
where, in accordance with the Lemma, $T_\iota=\calF\phi_\iota\in
\calS'(\bbR^n;M_m)$, and $Q(\partial/{\partial x})$ acts on $f_\iota$ in the sense of distributions. Indeed, if \eqref{eq9.10} and \eqref{eq9.11} hold, then
$$
\phi_\iota(\xi)=\frac1{(2\pi)^n}[Q(-i\xi)\widehat f_\iota(-\xi)+\widehat g_\iota(-\xi)],
$$
where $\widehat f_\iota$, $\widehat g_\iota$ are continuous and bounded on $\bbR^n$, and
$$
\sup\{\|\widehat f_\iota(\xi)\|_{M_m},\|\widehat g_\iota(\xi)\|_{M_m}
:\xi\in\bbR^n,\,\iota\in J\}<\infty,
$$
so that \eqref{eq9.9} is satisfied.
\goodbreak

{\it Construction of} $Q$, $f_\iota$ and $g_\iota$. We will construct $Q$, $f_\iota$ and $g_\iota$ in the form 
\begin{equation}\label{eq9.12}
Q=\Delta^k,\quad 2k\ge l+n+2,\ \quad f_\iota=T_\iota* u\quad\mbox{and}\ \quad
g_\iota=T_\iota*v 
\end{equation}
where $u\in C_K^l$, $v \in C_K^\infty$, $K=\{x\in\bbR^n:|x|\le 1\}$, and  $l$  is sufficiently large. Since $T_\iota*\varphi=\calF(\phi_\iota\bullet\calF^{-1}\varphi)$
for every $\varphi\in(\calS(\bbR^n))^m$, from \eqref{eq9.8} it follows that the family of convolution operators
$\{T_\iota\,*:\iota\in J\}\subset L((\calS(\bbR^n))^m)$ is equicontinuous. This implies that the family of convolution operators 
$$
\{(T_\iota\,*)|_{C_K^\infty}:\iota\in J\}\subset L(C_K^\infty;L^1(\bbR^n;M_m))
$$
is also equicontinuous. Since $L^1(\bbR^n;M_m)$ is a Banach space, there are $l\in\bbN_0$ and $C\in\]0,\infty\[$
such that
$$
\|T_\iota* w\|_{L^1(\bbR^n;M_m)}\le C\|w\|_{C_K^l}\ \quad\mbox{for every }
\iota\in J\mbox{ and }w\in C_K^\infty.
$$
Since $C_K^\infty$ is dense in $C_K^l$, it follows that for every $\iota\in J$
 the convolution $T_\iota* w$ of the distribution $T_\iota\in\calS(\bbR^n;M_m)$ 
and the scalar distribution represented by a function $w\in C_K^l$ is a function
$T_\iota* w\in L^1(\bbR^n;M_m)$ such that 
\begin{equation}\label{eq9.13}
\|T_\iota*w\|_{L^1(\bbR^n;M_m)}\le C\|w\|_{C_K^l}
\end{equation}	
where $C\in\]0,\infty\[$
is independent of $\iota\in J$ and $w\in C_K^l$.

{\it The formulas determining $ u$  and  $v$.}
 From \eqref{eq9.13} it follows that if $u\in C_K^l$
 and $v \in C_K^l$ are fixed, then for $f_\iota$ and $g_\iota$ defined by \eqref{eq9.12} the condition \eqref{eq9.10} is satisfied. Therefore the only problem which remains to be solved consists in choosing $u\in C_K^l$ and $v \in C_K^l$ so that 
$$
T_\iota=\Delta^k f_\iota+g_\iota=\Delta^k(T_\iota*u)+T_\iota*v=T_\iota*
(\Delta^k u+v).
$$
 To this end it is sufficient to choose $u\in C_K^l$, $v\in C_K^l$ and $k\in\bbN$ so that
$$
\Delta^k u+v=\delta.
$$

From the formula (II, 3;19) on p.~47 of \cite{S} (or from Theorem 5.1 of \cite{Ch}, p.~99, or from \cite{G-S1}, the example at the end of Sec.~III.2.1 of the 1958 ed. or Sec.~I.6.1 of the 1959 ed.) it follows that if $k\in\bbN$ and $2k\ge l+n+1$, then the differential operator $\Delta^k $ on $\bbR^n$ has a fundamental solution  $\calE$  which is a function belonging to $C^l(\bbR^n)$ such that $\calE|_{\bbR^n\setminus\{0\}}
\in C^\infty(\bbR^n\setminus\{0\})$. Let $\gamma\in C^\infty_K$ be equal to $1$ in a neighbourhood of~0. Define 
$$
u=\gamma\calE,\ \quad v=\Delta^k(1-\gamma)\calE.
$$
Then $u\in C_K^l$, $v\in C_K^\infty$ and $\Delta^k u+v=\Delta^k(\gamma\calE+(1-\gamma)\calE)
=\Delta^k\calE=\delta$.

\section{Proof of Theorem 1(i)}\label{sec10}

\noindent{\bf Proof of ${\rm(a)\Rightarrow(b)}\wedge(\omega_E\le\omega_0)$.} 
Let $E=(\calS(\bbR^n))^m$. Suppose that (a) holds, i.e. $\omega_0<\infty$. Then, by the implications \eqref{eq8.1}$\Rightarrow$\eqref{eq8.3} and \eqref{eq9.2}$\Rightarrow$\eqref{eq9.3}, for every $\varepsilon>0$ the family of multiplication operators
$$
\{e^{-(\omega+\varepsilon)t}e^{t\wid P}:t\in[0,\infty\[\}\subset L(E)
$$
is equicontinuous. Let $(G_t)_{t\in\bbR}\subset L((Z'_n)^m)$ be the one-parameter group \eqref{eq2.1}. By invariance of  $E$  with respect to the Fourier transformation it follows that
$S_t:=G_t|_E=\calF e^{t\wid P}\calF^{-1}\subset E$
for every $t\in[0,\infty\[$, and the family of operators
$(S_t)_{t\ge0}\subset L(E)$ is a one-parameter group for which $\omega_E\le\omega_0$. 
The Lemma from Section~\ref{sec7} shows that this
is a $(C_0)$-semigroup 
with infinitesimal generator $P(\partial/{\partial x})|_{(\calS(\bbR^n))^m}$. 

\medskip
\noindent{\bf Proof of ${\rm{(b)\Rightarrow}(a)}$  and 
${\rm(b)}\wedge(\omega_E<\infty)\Rightarrow(\omega_0\le\omega_E)$.} Suppose that (b) holds, i.e. $P(\partial/{\partial x})|_{(\calS(\bbR^n))^m}$ is the infinitesimal generator of a $(C_0)$-semigroup $(S_t)_{t\ge0}\subset L((\calS(\bbR^n))^m)$. Then, by the Uniqueness Theorem from Section~\ref{sec2}, $S_t=G_t|_{(\calS(\bbR^n))^m}$ for every
$t\in[0,\infty\[$ where the operators $G_t=\calF e^{t\wid P}\calF^{-1}\in L((Z'_n)^m)$,
$t\in\bbR$, constitute the one-parameter $(C_0)$-group 
$(G_t)_{t\in\bbR}\subset L((Z'_n)^m)$
considered in Section~\ref{sec2}. Consequently, if (b) holds, then
\begin{equation}\label{eq10.1}
e^{t\wid P}|_{(\calS(\bbR^n))^m}=\calF^{-1} S_t \calF |_{(\calS(\bbR^n))^m}\in L((\calS(\bbR^n))^m)
\end{equation}
for every $t\in[0,\infty\[$. By Theorem A of Section~\ref{sec9} it follows that
$$
e^{t\wid P}\in\calM((\calS(\bbR^n))^m)=\calO_M(\bbR^n;M_m)\ \quad
\mbox{for every }t\in[0,\infty\[.
$$
In particular,  $e^{\wid P}\in\calO_M(\bbR^n;M_m)$, and hence there is $k\in\bbN$ such that
$$
\sup\{(1+|\xi|)^{-k}\|e^{\wid P(\xi)}\|_{L(\bbC^m)}:\xi\in\bbR^n\}=K<\infty.
$$
By Corollary~2.4 on p.~252 of \cite{E-N}, this implies that
$$\eqalign{
\max\{\hRe \lambda:\lambda\in\sigma(\wid P(\xi))\}\le\log \rho(e^{\wid P(\xi)})
&\le\log\|e^{\wid P(\xi)}\|_{L(\bbC^m)}\cr
&\le\log K+k\log(1+|\xi|)\cr}
$$
for every $\xi\in\bbR^n$ where $\rho$ stands for the spectral radius. By the Corollary at the end of Section~\ref{sec5}, it follows that $\omega_0=\sup\{\hRe \lambda:\lambda\in\sigma(\wid P(\xi)),\,\xi\in\bbR^n\}\allowbreak<\infty$.

If (b) holds and $\omega_E<\infty$, then, in addition to \eqref{eq10.1}, for every
$\varepsilon>0$ the family of multiplication operators
$$\eqalign{
\{e^{-(\omega_E+\varepsilon)t} e^{t\wid P}|_{(\calS(\bbR^n))^m}:t\in[0,\infty\[\}
&=\{\calF^{-1}e^{-(\omega_E+\varepsilon)t}S_t\calF|_{(\calS(\bbR^n))^m}:t\in[0,\infty\[\}\cr
&\subset L((\calS(\bbR^n))^m)\cr}
$$
is equicontinuous, and hence, by the implication \eqref{eq9.3}$\Rightarrow$\eqref{eq9.2}, the condition \eqref{eq8.3} is satisfied for $\omega=\omega_E$. By the implication \eqref{eq8.3}$\Rightarrow$\eqref{eq8.1}, it follows that $\omega_0\le\omega_E$.

\section{Proof of Theorem 1(ii)}\label{sec11} 

The proof of Theorem 1(ii) is analogous to that of Theorem 1(i). 
The Lemma from Section~7 applies to $E=(\calS'_b(\bbR^n))^m=(\calF\calS'_b(\bbR^n))^m$
because $(Z_n)^m=(\calF C_c^\infty(\bbR^n))^m$ and $C_c^\infty(\bbR^n)$ is
sequentially dense in $\calS'_b(\bbR^n)$. This last may be proved by approximation of
distributions in $\calS'(\bbR^n)$ by cutting and regularizing. See \cite{R}, p.~253,
Proposition~4; \cite{T}, Sec.~28.

\section{Proof of Theorem 1(iii)}\label{sec12.}

The topology in $(C_b^\infty(\bbR^n))^m$ is determined by the sequence of norms
\begin{multline}\label{eq12.1}
\|u\|_j=\sup\{\|(\partial/\partial x)^\alpha u(x)\|_{\bbC^m}:\alpha\in\bbN_0^n,\,
|\alpha|\le j,\,x\in\bbR^n\},\\
j=0,1,\ldots,
\end{multline}  
where $u\in(C_b^\infty(\bbR^n))^m$. The implication ${\rm(a)\Rightarrow(b)}\wedge
(\omega_E\le\omega_0)$ for $E=(C_b^\infty(\bbR^n))^m$ will be proved by using some variants of estimates proved by I.~G.~Petrovski\u\i. These variants are uniform in $t\in[0,\infty\[$. The non-uniform 
estimates used in \cite{P} permit one only to prove
${\rm(a){\Rightarrow}(b)}$ without showing that $\omega_E<\infty$ and $\omega_E\le \omega_0$.

Suppose that $\omega_0<\infty$ and let $S_t^\circ:=G_t|_{(\calS(\bbR^n))^m} $  for every
$t\in[0,\infty\[$ where $(G_t)_{t\in\bbR}\subset L((Z_n')^m)$
is the one-parameter group \eqref{eq2.1}. Then, by Theorem~1(i), 
$(S_t^\circ)_{t\ge0}\subset L((\calS(\bbR^n))^m)$ is a one-parameter $(C_0)$-semigroup with infinitesimal generator $P(\partial/\partial x)|_{(\calS(\bbR^n))^m}$.

\medskip
\noindent
{\bf Lemma.} {\it If $\omega_0<\infty$, then for every $\varepsilon>0$ there is
$M_\varepsilon\in\]0,\infty\[$ such that
\begin{equation}\label{eq12.2}
\|(S_t^\circ u)(x)\|_{\bbC^m}\le M_{\varepsilon}\|u\|_{k_0}e^{(\omega_0+\varepsilon)t}
\prod_{\nu =1}^n(1+|x_\nu |)^{-2}
\end {equation} 
for every $u\in(C_I^\infty(\bbR^n))^m$, $t\in[0,\infty\[$ and $x=(x_1,\ldots,x_n)\in\bbR^n $ where $I=[-1/2,1/2]^n$ and $k_0=n((dm-1)(2n+1)+2)$.
}

\begin{pf}We follow I.~G.~Petrovski\u\i\  \cite{P}, pp.~12--17, but instead of analogues of \eqref{eq8.4} and \eqref{eq8.6} we use \eqref{eq8.1} and \eqref{eq8.3}. Notice that
\begin{equation}\label{eq12.3}
\min(1,a^{-k})\le 2^k(1+a)^{-k}\ \quad\mbox{for every }
\alpha\in\]0,\infty\[\mbox{ and }k\in\bbN.
\end {equation} 
Whenever $u\in(C_c^\infty(\bbR^n))^m$, $t\in[0,\infty\[$, $x=(x_1,\ldots,x_n)\in\bbR^n$ and 
$x_\nu \ne0$ for $\nu =1,\ldots,n$, then $\calF^{-1}u\in(\calS(\bbR^n))^m$, $e^{t\wid P}
\in\calO_M(\bbR^n;M_m)$, $e^{t\wid P}\calF^{-1}u\in(\calS(\bbR^n))^m$, and hence \vg
\begin{align}
(S_t^\circ u)(x)={}&(\calF e^{t\wid P}\calF^{-1}u)(x)
=\int_{\bbR^n} e^{i(x,\xi)} e^{t\wid P(\xi)}(\calF^{-1}u)(\xi)\,d\xi\notag\\
={}&(-1)^n(x_1\cdots x_n)^{-2}\notag\\
\label{eq12.4} 
&{}\times\int_{\bbR^n} e^{i(x,\xi)} \bigg(\frac{\partial^n}{\partial\xi_1
\cdots\partial\xi_n}\bigg)^2 [e^{t\wid P(\xi)}(\calF^{-1}u) (\xi)]\,d\xi.
\end{align}
By the Remark at the end of Section~\ref{sec8}, the inequality $\omega_0<\infty$ is equivalent to \eqref{eq8.3} with $\omega=\omega_0$ and $k_\alpha=(md-1)(|\alpha|+1)$. Therefore \eqref{eq12.3} and \eqref{eq12.4} imply that for every $\varepsilon>0$ there is $K_\varepsilon\in\]0,\infty\[$ such that
$$\displaylines{
\|(S_t^\circ u)(x)\|_{\bbC^m}
\le K_\varepsilon e^{(\omega_0+\varepsilon)t}\prod_{\nu =1}^n(1+|x_\nu |)^{-2}\hfill\cr
\hfill\times
\int_{\bbR^n}\prod_{\nu =1}^n(1+|\xi_\nu |)^{(md-1)(2n+1)}
\sup_{|\alpha|\le 2n}\|(\partial/\partial \xi)^\alpha(\calF^{-1}u)(\xi)\|_{\bbC^m}\,d\xi
\indent
\refstepcounter{equation}(\theequation)\label{eq12.5}
\cr}
$$
whenever $u\in(C_c^\infty(\bbR^n))^m$, $t\in[0,\infty\[$ and $x=(x_1,\ldots,x_n)\in\bbR^n$. Furthermore,
\begin{align*}
\bigg(\frac{\partial}{\partial\xi}\bigg)^\alpha(\calF^{-1}u)(\xi)
&=\frac1{(2\pi)^n}(-i)^{|\alpha|}\int_{\bbR^n}e^{-i\<x,\xi>}x^\alpha u(x)\,dx\\
&=\frac1{(2\pi)^n}(-i)^{|\alpha|}(\xi^\beta)^{-1}\int_{\bbR^n}e^{-i\<x,\xi>}
\bigg(\frac{\partial}{\partial x}\bigg)^\beta(x^\alpha u(x))\,dx
\end{align*}    
whenever $u\in(C_c^\infty(\bbR^n))^m$, $\alpha,\beta\in\bbN_0^n$, $\xi=(\xi_1,\ldots,\xi_n)\in\bbR^n$ and $\xi_\nu \ne 0$ for $\nu =1,\ldots, n$. Consequently, by \eqref{eq12.3}, for every
$k\in\bbN$ there is $C_k\in\]0,\infty\[$ such that whenever $u\in(C_I^\infty(\bbR^n))^m$, $\xi\in\bbR^n$ and $\alpha\in\bbN_0^n$, then
$$\displaylines{
\bigg\|\bigg(\frac{\partial}{\partial \xi}\bigg)^\alpha(\calF^{-1} u)(\xi)\bigg\|_{\bbC^m}\hfill\cr
\le C_k\prod_{\nu =1}^n(1+|\xi_\nu |)^{-k}\sup\{\|({\partial}/{\partial x})^\gamma
u(x)\|_{\bbC^m}:\gamma\in\bbN_0^n,\,|\gamma|\le kn,\,x\in I\}.
\hfill
\refstepcounter{equation}(\theequation)\label{eq12.6}
\cr}
$$
From \eqref{eq12.5} and \eqref{eq12.6} with $k=(md-1)(2n+1)+2$ one obtains \eqref{eq12.2}.   \end{pf}

\noindent{\bf Proof of ${\rm(a)\Rightarrow(b)}\wedge(\omega_E\le\omega_0)$  for $E=(C_b^\infty(\bbR^n))^m$.} Let $\bbZ$ be the set of all integers (positive or non-positive). 
For any $z=(z_1,\ldots,z_n)\in\bbZ^n$ denote by $\tau_z$ the operator of translation:
$(\tau_z f)(x):=f(x_1+\frac12z_1,\ldots,x_n+\frac12z_n)$ 
for every function  $f$ defined on $\bbR^n$ and every $x=(x_1,\ldots,x_n)\in\bbR^n$. Let 
$I=[-1/2,1/2]^n$. Following \cite{P}, fix a function $v \in C_I^\infty(\bbR^n)$ with values in $[0,1]$
such that $\sum_{z\in\bbZ^n}\tau_zv \equiv 1$  on $\bbR^n$. Assume that 
$\omega_0<\infty$. Since the operators $S_t^\circ$, $\tau_z$ and $(\partial/\partial x)^\alpha$ commute, and
$$ 
\frac \partial{\partial t}S_t^\circ(u\tau_z v )=P\bigg(\frac \partial{\partial x}\bigg)
S_t^\circ(u\tau_z v ),
$$
one has
\begin{equation}\label{eq12.7}
\bigg(\frac \partial{\partial x}\bigg)^\alpha\bigg(\frac \partial{\partial t}\bigg)^\beta
S_t^\circ(u\tau_z v )=\tau_zS_t^\circ
\bigg(\frac \partial{\partial x}\bigg)^\alpha\bigg(P\bigg(\frac \partial{\partial x}\bigg)\bigg)^\beta(v \tau_{-z} u)
\end{equation}
for every $u\in(C_b^\infty(\bbR^n))^m$, $t\in[0,\infty\[$, $\alpha\in\bbN_0^n$
 and $\beta\in \bbN_0$. Since the norms \eqref{eq12.1} are translation-invariant, from \eqref{eq12.2} and \eqref{eq12.7} it follows that 
$$\displaylines{\indent
\mbox{whenever }u\in(C_b^\infty(\bbR^n))^m,\ t\in[0,\infty\[,\ x=(x_1,\ldots,x_n)\in\bbR^n,\hfill\cr 
\indent z=(z_1,\ldots,z_n)\in\bbZ^n,\  \alpha\in\bbN_0^n,\ \beta\in\bbN_0\mbox{ and }
\varepsilon>0,\mbox{ then}\hfill\cr 
\indent\bigg\|\bigg[\bigg(\frac \partial{\partial x}\bigg)^\alpha\bigg(\frac \partial{\partial t}\bigg)^\beta S_t^\circ(u\tau_z v )\bigg](x)\bigg\|_{\bbC^m}\hfill\cr 
\hfill
\le M_{\varepsilon,\alpha,\beta}\|u\|_{k_0+|\alpha|+\beta d}
e^{(\omega_0+\varepsilon)t}\prod_{\nu =1}^n(1+|x_\nu +\tfrac12 z_\nu |)^{-2}\indent
\refstepcounter{equation}(\theequation)\label{eq12.8}\cr
\indent\mbox{where $k_0=n((dm-1)(2n+1)+2)$ and $M_{\varepsilon,\alpha,\beta}\in\]0,\infty\[$ depends}\hfill\cr
\indent\mbox{only on  $\varepsilon$, $\alpha$ and $\beta$.}\hfill\cr}
$$
The series $\sum_{z\in\bbZ^n}\prod_{\nu =1}^n(1+|x_\nu +\tfrac12 z_\nu |)^{-2}$ of continuous functions  of $x=(x_1,\ldots,x_n)$ is uniformly convergent on every bounded subset 
of $\bbR^n$, so that its sum
$$
s(x)=\sum_{z\in\bbZ^n}\prod_{\nu =1}^n(1+|x_\nu +\tfrac12 z_\nu |)^{-2}
$$
\looseness2
is a continuous function of~$x$. Since  $s$  is periodic, it is bounded on $\bbR^n$. Consequently, if $u\in(C_b^\infty(\bbR^n))^m$
and $u_z(t,x)=[S_t^\circ(u\tau_zv )](x)$, then, by \eqref{eq12.8}, 
 $\sum_{z\in\bbZ^n} u_z(t,x)$ 
is a series of functions of $(t,\!x)$ belonging to $C^\infty([0,\!\infty\[;\allowbreak
(\calS(\bbR^n))^m)$, and it converges uniformly on every bounded subset of
$[0,\infty\[\times\bbR^n$ together with all partial derivatives in  $t$  and $x_1,\ldots,x_n$. Therefore, by the theorem on term by term differentiation, whenever $u\in(C_b^\infty(\bbR^n))^m$ and 
\begin{equation}\label{eq12.9}
u(t,x):=\sum_{z\in\bbZ^n} u_z(t,x)\ \quad\mbox{for }(t,x)\in[0,\infty\[\times\bbR^n,
\end{equation}
then $u(\bullet,\bullet)\in (C^\infty([0,\infty\[\times\bbR^n))^m$, $u(\bullet,\bullet)$   satisfies the PDE 
$$
\frac \partial{\partial t}u(t,x)=P\bigg(\frac \partial{\partial x}\bigg)u(t,x)\ \quad
\mbox{for }(t,x)\in[0,\infty\[\times\bbR^n,
$$
\looseness2
 and $u(0,x)=u(x)$ for $x\in\bbR^n$. Inequality \eqref{eq12.8} and boundedness of  $s$  imply that\vg 
$$\displaylines{\indent
\mbox{for every $\varepsilon>0$, $j\in\bbN_0$, 
and $\beta\in\bbN_0$ there is}\hfill\cr
\indent M_{\varepsilon,j,\beta}\in\]0,\infty\[
\mbox{ such that}\hfill\cr
\indent
\bigg\|
\bigg(\frac \partial{\partial t}\bigg)^\beta u(t,\bullet)\bigg\|_{j}\le 
M_{\varepsilon,j,\beta}\|u\|_{j+\beta d+k_0}e^{(\omega_0+\varepsilon)t}\hfill
\refstepcounter{equation}(\theequation)\label{eq12.10}\cr
\indent\mbox{whenever $u\in(C_b^\infty(\bbR^n))^m$ and $t\in[0,\infty\[$.}\hfill
\cr}
$$
Consequently,
$$\displaylines{\indent
\mbox{if } u\in(C_b^\infty(\bbR^n))^m,\mbox{ then }u(\bullet,\bullet)\in C^\infty([0,\infty\[;(C_b^\infty(\bbR^n))^m)\hfill 
\refstepcounter{equation}(\theequation)\label{eq12.11}
\cr} 
$$                  
and
$$\displaylines{\indent
\mbox{if  $u\in(C_b^\infty(\bbR^n))^m$ and $t_0\in[0,\infty\[$, 
 then}\hfill\cr
\indent	 \lim_{[0,\infty\[
\ni t\to t_0}\frac 1{t-t_0}(u(t,\bullet)-u(t_0,\bullet))
=P\bigg(\frac{\partial}{\partial x}\bigg) u (t_0,\bullet)\hfill
\refstepcounter{equation}(\theequation)\label{eq12.12}\cr
\indent\mbox{in  the topology of $(C_b^\infty(\bbR^n))^m$.}\hfill\cr}
$$                                 

{\spaceskip.33em plus.22em minus.17em 
Let $(G_t)_{t\in\bbR}\!\subset\! L((Z_n')^m)$ be the one-parameter group \eqref{eq2.1}. 
By the Uniqueness Theorem from Section~2, from \eqref{eq12.11}, \eqref{eq12.12} and the continuous imbeddings $(C_b^\infty(\bbR^n))^m
\subset(\calS'_b(\bbR^n))^m\subset(Z'_n)^m$
it follows that
\begin{equation}\label{eq12.13}
u(t,\bullet)=G_tu\ \quad\mbox{for every }
t\in[0,\infty\[\mbox{ and }u\in(C_b^\infty(\bbR^n))^m.
\end{equation}
Consequently, by \eqref{eq12.10}, one has $G_t|_{C_b^\infty(\bbR^n)^m}
\in L((C_b^\infty(\bbR^n))^m)$, and, by \eqref{eq12.11}--\eqref{eq12.13},
$(S_t)_{t\ge 0}:=(G_t|_{(C_b^\infty(\bbR^n))^m})_{t\ge 0}
\subset L((C_b^\infty(\bbR^n)^m)$ is a $(C_0)$-semigroup 
with infinitesimal generator $P(\partial/\partial x)|_{(C_b^\infty(\bbR^n))^m}$.

\medskip

\noindent
{\bf Proof of ${\rm(a)\Rightarrow(b)}$ and ${\rm(b)}\wedge
(\omega_E<\infty)\Rightarrow
(\omega_0\le\omega_E)$ for
$E=(C_b^\infty(\bbR^n))^m$.}
 We adapt an argument due to I.~G.~Petrovski\u\i\  in \cite{P}, pp.~7--9, to the semi\-group-theoretical formulation. For every $\xi\in\bbR^n$ let $\chi_\xi$ be the character of
$\bbR^n$  such that $\chi_\xi(x)=e^{i(x,\xi)}$ for $x\in\bbR^n$. The Fourier transformation  $\calF$  is an isomorphism of $\calD'(\bbR^n)$ onto $Z'_n$, and since it acts on 
$(\calD'(\bbR^n))^m$ coordinatewise, it  is also an isomorphism of
$(\calD'(\bbR^n))^m$ onto $(Z'_n)^m$. One has
$\chi_\xi\in C_b^\infty(\bbR^n)\subset \calS'(\bbR^n)\subset Z'_n$ and
$\chi_\xi=\calF\delta_\xi$. Whenever $t\in\bbR$ is fixed, then
$e^{t\wid P}\in C^\infty(\bbR^n;M_m)$ and the multiplication operator
$e^{t\wid P}\bullet$ maps $(\calD'(\bbR^n))^m$ into $(\calD'(\bbR^n))^m$. Consequently, whenever $t\in\bbR^n$, $\xi\in\bbR^n$ and $z\in\bbC^m$, then 
$$\eqalign{
[\calF(e^{t\wid P}\bullet)\calF^{-1}](\chi_\xi\otimes z)
&=\calF(e^{t\wid P}\bullet(\delta_\xi\otimes z))\cr
&=\calF(\delta_\xi\otimes (e^{t\wid P(\xi)}z))
=\chi_\xi\otimes(e^{t\wid P(\xi)}z).\cr}
$$
This means that
\begin{equation}\label{eq12.14}	 
G_t(\chi_\xi\otimes z)=\chi_\xi\otimes(e^{t\wid P(\xi)}z)\ \quad\mbox{for }
t\in\bbR^n,\,\xi\in\bbR^n\mbox{ and }z\in\bbC^m
\end{equation}
where $(G_t)_{t\in\bbR}\subset L((Z'_n)^m)$ is the one-parameter $(C_0)$-group \eqref{eq2.1}.

Suppose now that (b) holds for $E=(C_b^\infty(\bbR^n))^m$, and let
$(S_t)_{t\ge0}\subset\break L((C_b^\infty(\bbR^n))^m)$ be the one-parameter $(C_0)$-semigroup with infinitesimal generator $P(\partial/\partial x)|_{(C_b^\infty(\bbR^n))^m}$. 
Then, by the Uniqueness Theorem from Section~\ref{sec2},
$G_t|_{(C_b^\infty(\bbR^n))^m}=S_t$ for every $t\in[0,\infty\[$. Consequently if
$t\in[0,\infty\[$ and $u\in(C_b^\infty(\bbR^n))^m$, then $G_tu\in(C_b^\infty(\bbR^n))^m$, 
and for every $t\in[0,\infty\[$ there are $K_t\in\]0,\infty\[$ and $k_t\in\bbN_0$ such that
\begin{equation}\label{eq12.15}
\|G_tu\|_0\le K_t\|u\|_{k_t}\quad\mbox{whenever}\quad
u\in(C_b^\infty(\bbR^n))^m.
\end{equation}  
From \eqref{eq12.14} and \eqref{eq12.15} it follows that if $t\in[0,\infty\[$ and 
$\xi\in\bbR^n$, then\vg
\begin{align*}
 \|e^{t\wid P(\xi)}\|_{L(\bbC^m)}&=\sup\{\|\chi_\xi(x)\otimes( e^{t\wid P(\xi)}z)\|_{\bbC^m}
:x\in\bbR^n,\,z\in \bbC^m,\,\|z\|_{\bbC^m} \le1\}\\
&= \sup_{\|z\|_{\bbC^m} \le1} \|\chi_\xi\otimes( e^{t\wid P(\xi)}z)\|_{0}
= \sup_{\|z\|_{\bbC^m} \le1}\|G_t (\chi_\xi\otimes z)\|\\
&\le \sup_{\|z\|_{\bbC^m}\le1} K_t\|\chi_\xi\otimes z\|_{k_t}\\
&\le K_t\sup\{|(\partial/\partial x)^{\alpha}\chi_\xi(x)|:x\in\bbR^n,\,\alpha\in\bbN_0^n,\,
|\alpha|\le k_t\}\\       
&\le K_t\sup\{|\xi_1^{\alpha_1}\cdots \xi_n^{\alpha_n}|:\alpha_\nu \in\bbN_0\mbox{ for }\nu =1,\ldots,n, \\
&\hskip220pt \alpha_1+\cdots+ \alpha_n\le k_t\}\\
 &\le K_t(1+|\xi_1|+\cdots+ |\xi_n|)^{k_t}\le K_t(1+n |\xi|)^{k_t}
\le  K_t(1+|\xi| )^{nk_t}.       
\end{align*}                 
By Corollary 2.4, p.~252 of \cite{E-N} it follows that whenever $t\in\]0,\infty\[$ and
$\xi\in\bbR^n$, then
\begin{align}
\max\{\hRe\lambda:\lambda\in\sigma(\wid P(\xi))\}
&= t^{-1}\log\rho(e^{t\wid P(\xi)})
\le t^{-1}\log\|e^{t\wid P(\xi)}\|_{L(\bbC^m)}\notag\\
\label{eq12.16}
&\le t^{-1}\log K_t+ t^{-1}nk_t\log(1+|\xi|) 
\end{align}
where $\rho(e^{t\wid P(\xi)})$ denotes the spectral radius of the matrix $e^{t\wid P(\xi)}$. From \eqref{eq12.16}, by the Corollary at the end of Section~\ref{sec5}, it follows that  $\omega_0<\infty$, where
$$
\omega_0:=\sup\{\hRe\lambda:\lambda\in\sigma(\wid P(\xi)),\,\xi\in\bbR^n\}<\infty,
$$
proving (a).

If ${\rm(b)}\wedge(\omega_E<\infty)$ holds for $E=(C_b^\infty(\bbR^n))^m$, 
then there is a $k\in\bbN$, and for every $\varepsilon>0$ there is an
$M_\varepsilon\in\]0,\infty\[$ such that for $k_t$ and $K_t$ occurring in \eqref{eq12.15} and \eqref{eq12.16} one has
$$ 
k_t\le k\quad\mbox{and}\quad K_t\le M_\varepsilon e^{(\omega_E+\varepsilon)t}
\ \quad\mbox{for every }t\in[0,\infty\[.
$$
Consequently, by \eqref{eq12.16},
$$ 
\max\{\hRe\lambda:\lambda\in\sigma(\wid P(\xi))\}
\le\omega_E+\varepsilon+t^{-1}[\log M_\varepsilon+nk\log(1+|\xi|)]
$$ 
for every $\xi\in\bbR^n$ and $t>0$. Since $\varepsilon>0$ and
$t>0$ are arbitrary, it follows that 
$\max\{\hRe\lambda:\lambda\in\sigma(\wid P(\xi))\}
\le\omega_E$ for every $\xi\in\bbR^n$, proving that $\omega_0\le\omega_E$.

\end{document}